\documentclass[12pt,a4paper,twoside]{article}
\usepackage[english]{babel}
\usepackage{latexsym,amsfonts,amssymb,amsmath,longtable,amsthm}
\renewcommand{\le}{\leqslant}
\renewcommand{\ge}{\geqslant}

\begin{document}

\newtheorem{ttt}{Theorem}[section]
\newtheorem{lem}[ttt]{Lemma}
\newtheorem{cor}[ttt]{Corollary}
\newtheorem{prop}[ttt]{Proposition}

\theoremstyle{definition}
\newtheorem{df}[ttt]{Definition}
\newtheorem{ex}[ttt]{Example}
\newtheorem{xca}[ttt]{Excersize}
\newtheorem{conj}[ttt]{Conjecture}
\newtheorem{tab}{Table}

\theoremstyle{remark}
\newtheorem{rem}[ttt]{Remark}
\newtheorem{problem}[ttt]{Problem}
\renewcommand{\mod}{\mathrm{mod\ }}

\begin{center}
%\underline{Sobolev Institute of Mathematics\hspace{4.2cm}Preprint
%\textcyrillic{\No}~152, March 2005} \vspace{1cm}

{\bf \Large An adjacency criterion for the prime graph of a~finite
simple group}\footnote{The authors were supported by the~Russian
Foundation for Basic Research (Grant 05--01--00797),
the State Maintenance Program for the Leading Scientific Schools of
the Russian Federation (Grant NSh--2069.2003.1), the Program
``Development of the Scientific Potential of Higher School'' of the
Ministry for Education of the Russian Federation (Grant~8294),
the Program ``Universities of Russia'' (Grant UR.04.01.202), and by Presidium
SB RAS Grant 86-197.}
\vspace{1cm}

{\sc A. V. Vasil'ev, E. P. Vdovin}
\end{center}
\vspace{1cm}

{\small{\bfseries Abstracts:} In the paper we give an exhaustive
arithmetic criterion of adjacency in prime graph of every
finite nonabelian simple group. By using this criterion for the prime graph of
every finite simple groups an independence set with the maximal number of
vertices, an independence set containing $2$ with the maximal number
of vertices, and the orders of these independence sets are given. We
assemble this information in the tables at the end of the paper.
Several applications of obtained results for various problems of
finite group theory are considered, in particular, on the recognition problem
of a finite group by spectrum.\vspace{0.2cm}

{\bfseries Keywords:} finite groups, finite simple groups, groups of
Lie type, spectrum of finite group, recognition by spectrum, prime
graph of finite group, independent number of prime graph,
2-independent number of prime graph} \vspace{0.5cm}

Let $G$ be a finite group, $\pi(G)$ be the set of all prime divisors
of its order and $\omega(G)$ be the spectrum of $G$, that is the set
of all its element orders. A graph $GK(G)=\langle
V(GK(G)),E(GK(G))\rangle$, where $V(GK(G))$ is the set of vertices
and $E(GK(G))$ is the set of edges, is called the {\em
Gruenberg-Kegel graph} (or the {\em prime graph}) of $G$ if
$V(GK(G))=\pi(G)$ and the edge $(r,s)$ is in $E(GK(G))$ if and only
if $rs\in\omega(G)$. Primes $r,s\in \pi(G)$ are called {\em
adjacent} if they are adjacent as vertices of $GK(G)$, that is
$(r,s)\in E(GK(G))$. Otherwise $r$ and $s$ are called {\em
non-adjacent}.

The properties of the prime graph $GK(G)$ yield a rich information
on the structure of $G$ (see [1--4] and Sections \ref{influence} and
\ref{applications} of the paper). The main purpose of this article
is to give an exhaustive arithmetic criterion of adjacency in prime
graph $GK(G)$ for every finite nonabelian simple group $G$. Sections
\ref{prlmn}--\ref{two} are devoted to this goal. In
Section~\ref{influence} we discuss some recent results on prime
graph of finite groups that gave us the motivation for the present
work. Furthermore, we explain the importance of so-called
independence numbers of prime graph for investigations on a group
structure. In Section~\ref{sporadicalt} we calculate those
invariants for all finite nonabelian simple groups (the resulting
tables are assembled in Section~\ref{tables}). The applications of
our results are considered in Section \ref{applications}.

\section{Preliminaries}\label{prlmn}
\pagestyle{myheadings} \markboth{\underline{\hspace{5.1cm} A. V. Vasil'ev, E. P.
Vdovin\hspace{5.1cm}}}{\underline{\hspace{4.8cm}Prime graph of finite simple group\hspace{4.8cm}
\hspace{-5mm}}}

Our notations is standard.  If $n$ is a natural number, $\pi$ a set
of primes, then by $\pi(n)$ we denote the set of all prime divisors
of $n$, by $n_\pi$ we denote the maximal divisor $t$ of $n$ such
that $\pi(t)\subseteq\pi$. Note that for a finite group $G$,
$\pi(G)=\pi(|G|)$ by definition.

The adjacency criterion of two prime divisors  for alternating
groups is obvious and can be given as follows.

\begin{prop}\label{adjalt}
Let $G=Alt_n$ be an alternating group of degree $n$.

\begin{itemize}
\item[{\em (1)}] Let $r,s\in\pi(G)$ be odd primes. Then $r,s$
are non-adjacent if and only if $r+s>n$.
\item[{\em (2)}] Let $r\in\pi(G)$ be an odd prime. Then $2,r$ are
non-adjacent if and only if $r+4>n$.
\end{itemize}
\end{prop}

The information on the adjacency of vertices in a prime graph for
every sporadic group and Tits group ${}^2F_4(2)'$ can be extracted
from \cite{Atlas} or~\cite{GAP}. Thus, we need to consider only
simple groups of Lie type.

For groups of Lie type and linear algebraic groups our notations
agrees with that in \cite{Ca1} and \cite{Hu1} respectively. Denote
by $G_{sc}$ a universal group of Lie type. Then every factor group
$G_{sc}/Z$, where $Z\leq Z(G_{sc})$, we call a {\em group of Lie
type}. Almost in all cases $G_{sc}/Z(G_{sc})$ is simple and we say
that $G_{ad}=G_{sc}/Z(G_{sc})$ is of {\em adjoint type}. Some groups
of Lie type over small fields are not simple. In
Table~\ref{notsimple} we assemble the information of \cite[Theorems
11.1.2 and~14.4.1]{Ca1} and give all such exceptions.  Sometimes we
shall use notations $A_n^\varepsilon(q), D_n^\varepsilon(q)$, and
$E_6^\varepsilon(q)$, where $\varepsilon\in\{+,-\}$, and
$A_n^+(q)=A_n(q)$, $A_n^-(q)={^2A_n(q)}$, $D_n^+(q)=D_n(q)$,
$D_n^-(q)={^2D_n(q)}$, $E_6^+(q)=E_6(q)$, $E_6^-(q)={^2E_6(q)}$.

If $G$ is isomorphic to ${^2A_n(q)}$, ${^2D_n(q)}$, or ${^2E_6(q)}$
we say that $G$ is defined over $GF(q^2)$, if $G\simeq{^3D_4(q)}$ we
say that $G$ is defined over $GF(q^3)$ and we say that $G$ is
defined over $GF(q)$ for other finite groups of Lie type. The field
$GF(q)$ in all cases is called the {\em base field} of $G$. If $G$
is a universal finite group of Lie type with the base field $GF(q)$,
then there are a natural number $N$($=\vert\Phi^+\vert$ in most
cases) and a polynomial $f(t)\in\mathbb{Z}[t]$ such that
$|G|=f(q)\cdot q^N$ and $(q,f(q))=1$ (see \cite[Theorems~9.4.10
and~14.3.1]{Ca1}). This polynomial we denote by $f_G(t)$. If $G$ is
not universal then there is a universal group $K$ with $G=K/Z$, where $Z=Z(K)$,
and we define~${f_G(t)=f_K(t)}$.

Assume that $\overline{G}$ is a connected simple linear algebraic
group defined over an algebraically closed field of a positive
characteristic $p$. Let $\sigma$ be an endomorphism of
$\overline{G}$ such that
$\overline{G}_\sigma=C_{\overline{G}}(\sigma)$ is a finite set. Then
$\sigma$ is said to be a {\em Frobenius map} and
$G=O^{p'}(\overline{G}_\sigma)$ is appeared to be a finite group of
Lie type. Moreover all finite groups of Lie type both split and
twisted can be derived in this way. Below we assume that for every
finite group $G$ of Lie type we fix (in some way) a linear algebraic
group $\overline{G}$ and a Frobenius map $\sigma$ such that
$G=O^{p'}(\overline{G}_\sigma)$. If $\overline{G}$ is simply
connected, then $G=\overline{G}_\sigma=
O^{p'}(\overline{G}_\sigma)$; and if $\overline{G}$ is of adjoint
type, then $\overline{G}_\sigma=\widehat{G}$ is the group of
inner-diagonal automorphisms of~$G$ (see~{\cite[\S~12]{St}}).

If $\overline{R}$ is a $\sigma$-stable reductive subgroup of
$\overline{G}$, then $\overline{R}_\sigma\cap G=\overline{R}\cap G$
is called a {\em reductive subgroup} of $G$. If $\overline{R}$ is of
maximal rank, then $\overline{R}_\sigma\cap G$ is also said to be  of {\em
maximal rank}. Note that if $\overline{R}$ is $\sigma$-stable
reductive of maximal rank, then $\overline{R}=\overline{T}\cdot
\overline{G}_1\ast\ldots\ast\overline{G}_k$, where $\overline{T}$ is
some $\sigma$-stable maximal torus of $\overline{R}$ and
$\overline{G}_1,\ldots,\overline{G}_k$ are subsystem subgroups
of~$\overline{G}$. Furthermore
$\overline{G}_1\ast\ldots\ast\overline{G}_k=\left[\overline{R},
\overline{R}\right]$. It is known that
$\overline{R}_\sigma=\overline{T}_\sigma G_1\ast\ldots\ast G_m$,
subgroups $G_1,\ldots,G_m$ we call {\em subsystem subgroups} of~$G$.
In general $m\le k$ and for all $i$ the base field of $G_i$ is
equal to $GF(q^{\alpha_i})$, where $\alpha_i\ge1$. There is a nice
algorithm due to Borel and de Siebental determining all subsystem
subgroups of $\overline{G}$ (see \cite{BorSie} and also \cite{Dyn}).
One has to consider the extended Dynkin diagram of $\overline{G}$
and remove any number of nodes. Connected components of the
remaining graph are Dynkin diagrams of subsystem subgroups of
$\overline{G}$ and Dynkin diagrams of all subsystem subgroups can be
derived in this way.

If $\overline{T}$ is a $\sigma$-stable torus of $\overline{G}$ then
$T=\overline{T}\cap G=\overline{T}_\sigma\cap G$ is called a {\em
torus} of $G$. If $\overline{T}$ is maximal, then $T$ is a {\em
maximal torus} of $G$. If $G$ is neither a Suzuki group, nor a Ree
group, then for every maximal torus $T$ we have that $\vert
\overline{T}_\sigma\vert=g(q)$, where $GF(q)$ is the base field of
$G$, $g(t)$ is a polynomial of degree $n$ dividing $f_G(t)$ and $n$
is the rank of $\overline{G}$.  For more details see~\cite[Chapter~1]{Car6}.

\begin{tab}\label{notsimple}
{\bfseries Non simple groups of Lie type}
\begin{longtable}{|c|l|}\hline
Group&Properties\\ \hline
$A_1(2)$&Group is soluble\\
$A_1(3)$&Group is soluble\\
$B_2(2)$&$B_2(2)\simeq\mathrm{Sym}_6$\\
$G_2(2)$&$[G_2(2),G_2(2)]\simeq {}^2A_2(3)$\\
${}^2A_2(2)$&Group is soluble\\
${}^2B_2(2)$&Group is soluble\\
${}^2G_2(3)$&$[{}^2G_2(3),{}^2G_2(3)]\simeq A_1(8)$\\
${}^2F_4(2)$&$[{}^2F_4(2),{}^2F_4(2)]={}^2F_4(2)'$ is the Tits
group\\ \hline
\end{longtable}
\end{tab}

In Lemmas \ref{toriofclassgrps} and
\ref{toriofexcptgrps} we assemble the information about maximal tori in finite
 simple groups of Lie type.
% We assume that $G$ is simple, i.~e., we do not
%consider groups from Table~\ref{notsimple}.

\begin{lem}\label{toriofclassgrps} {\em (see \cite[Propositions~7--10]{Ca2} and \cite{Ca3})}
Let $\overline{G}$ be a connected simple classical algebraic group
of adjoint type and let $G=O^{p'}(\overline{G}_\sigma)$ be the finite
simple classical group.
\begin{itemize}
\item [{\em 1.}] Every
maximal torus $T$ of $G=A_{n-1}^\varepsilon(q)$ has the order
\begin{equation*}
\frac{1}{(n,q-(\varepsilon1))(q-(\varepsilon1))}
(q^{n_1}-(\varepsilon1)^{n_1})\cdot
(q^{n_2}-(\varepsilon1)^{n_2})\cdot\ldots\cdot
(q^{n_k}-(\varepsilon1)^{n_k})
\end{equation*}
for appropriate partition
$n_1+n_2+\ldots+n_k=n$ of~$n$. Moreover, for every partition there exists a
torus of corresponding order.

\item [{\em 2.}] Every maximal torus $T$ of $G$, where $G=B_n(q)$ or
$G=C_n(q)$, has the order
\begin{equation*}
\frac{1}{(2,q-1)}
(q^{n_1}-1)\cdot(q^{n_2}-1)\cdot\ldots\cdot(q^{n_k}-1)\cdot
(q^{l_1}+1)\cdot (q^{l_2}+1)\cdot\ldots\cdot(q^{l_m}+1)
\end{equation*}
for appropriate partition $n_1+n_2+\ldots+n_k+l_1+l_2+\ldots+l_m=n$ of~$n$.
Moreover, for every partition there exists a
torus of corresponding order.

\item [{\em 3.}] Every maximal torus $T$ of $G=D_n^\varepsilon(q)$ has the order
\begin{equation*}
\frac{1}{(4,q^n-\varepsilon1)}\cdot (q^{n_1}-1)\cdot
(q^{n_2}-1)\cdot \ldots\cdot
(q^{n_k}-1)\cdot(q^{l_1}+1)\cdot
(q^{l_2}+1) \cdot\ldots\cdot (q^{l_m}+1)
\end{equation*}
for appropriate partition
$n_1+n_2+\ldots+n_k+l_1+l_2+\ldots+l_m=n$ of $n$,
where $m$ is even if $\varepsilon=+$ and $m$ is odd if~$\varepsilon=-$.
Moreover, for every partition there exists a
torus of corresponding order.
\end{itemize}
\end{lem}

\begin{lem}\label{toriofexcptgrps} {\em (see \cite{Ca3} and \cite{der1})}
Let $\overline{G}$ be a connected simple exceptional algebraic group
of adjoint type and let $G=O^{p'}(\overline{G}_\sigma)$ be the finite
simple exceptional group of Lie type.
\begin{itemize}
\item[{\em 1.}]
Every maximal torus $T$ of $G=G_2(q)$ has one of the following orders:
\begin{equation*}
(q\pm1)^2, q^2-1, q^2\pm q+1.
\end{equation*}
Moreover, for every number given above there exists a torus of corresponding
order.
\item [{\em 2.}] Every maximal torus $T$ of $G=F_4(q)$ has one of the
following orders:
\begin{equation*}
(q\pm1)^4, (q\pm1)^2\cdot(q^2\pm1), (q^2\pm1)^2, (q\pm1)(q^3\pm1),
q^4\pm1, (q^2\pm q+1)^2, q^4-q^2+1.
\end{equation*}
Here and below symbol $\pm$ means that we can choose either ``$+$'' or
``$-$'' independently for all multiplies, i.~e., $(q\pm1)^2(q^2\pm1)$
is equal to either $(q-1)^2(q^2-1)$, or $(q+1)^2(q^2-1)$, or
$(q-1)^2(q^2+1)$, or $(q+1)^2(q^2+1)$. Moreover, for every number given above
there exists a torus of corresponding
order.

\item[{\em 3.}]
For every maximal torus $T$ of $G=E_6^\varepsilon(q)$, the number
$(3,q-\varepsilon1)\vert T\vert$ is equal to one of the following:

%\begin{multline*}
$(q-\varepsilon1)^k\cdot(q+\varepsilon1)^{6-k}\text{,
   }2\le k\le6;
  (q^k-(\varepsilon1)^k)\cdot
  (q^{6-k}-(\varepsilon1)^{6-k})\text{,  }1\le k \le5;$\par
$  (q^k-(\varepsilon1)^k)\cdot (q-\varepsilon1)^{6-k}\text{,
  }3\le k\le 6; (q^3-\varepsilon1)(q^2-1)(q\pm1);
  (q^5-\varepsilon1)(q+\varepsilon1);$\par
$  (q^3+\varepsilon1)(q^2\pm1)(q-\varepsilon1); (q^4+1)(q^2-1);
  (q^2+1)^2(q-\varepsilon1)^2; (q^2+\varepsilon q+1)^3;$\par
$  (q^2+\varepsilon q+1)^2(q^2-1); (q^4-1)(q+\varepsilon 1)^2;
  (q^3+\varepsilon1)(q^2+\varepsilon q+1)(q+\varepsilon 1);$\par
$  (q^4-q^2+1)(q^2+\varepsilon q+1); q^6+\varepsilon q^3+1;
  (q^2+\varepsilon q+1)(q^2-\varepsilon q+1)^2.$
%\end{multline*}
Moreover, for every number $n$ given above there exists a torus $T$ with
$(3,q-\varepsilon1)\vert T\vert=n$.
\item[{\em 4.}] For every maximal torus $T$ of $G=E_7(q)$, the number
  $m=(2,q-1)\vert T\vert$ is equal to one of the following:
$(q+1)^{n_1}(q-1)^{n_2},$ $n_1+n_2=7;$ $(q^2+1)^{n_1}(q+1)^{n_2}(q-1)^{n_3},$
$1\leqslant
n_1\leqslant2,$ $2n_1+n_2+n_3=7,$ and $m\neq(q^2+1)(q\pm1)^5;$
$(q^3+1)^{n_1}(q^3-1)^{n_2}(q^2+1)^{n_3}(q+1)^{n_4}(q-1)^{n_5},$ $1\leqslant
n_1+n_2\leqslant2,$
$3n_1+3n_2+2n_3+n_4+n_5=7,$ and $m\neq(q^3+\epsilon1)(q-\epsilon1)^4,$
$m\neq(q^3\pm1)(q^2+1)^2,$
$m\neq(q^3+\epsilon1)(q^2+1)(q+\epsilon1)^2;$ 
$(q^4+1)(q^2\pm1)(q\pm1);$ $(q^5\pm1)(q^2-1);$ $(q^5+\epsilon1)(q+\epsilon1)^2;$
$q^7\pm1;$
$(q-\epsilon1)\cdot (q^2+\epsilon q+1)^3; (q^5-\epsilon1)\cdot
  (q^2+\epsilon q +1); (q^3\pm1)\cdot (q^4-q^2+1); (q-\epsilon1)\cdot
  (q^6+\epsilon q^3+1);$ $(q^3-\epsilon1)\cdot(q^2-\epsilon
  q+1)^2,$ where $\epsilon=\pm$. Moreover, for every number $m$ given above
there exists a torus $T$ with $(2,q-1)\vert T\vert=m$.
$(2,q-1)\vert T\vert=n$.
\item[{\em 5.}] Every maximal torus $T$ of $G=E_8(q)$ has one of the
  following orders:
$(q+1)^{n_1}(q-1)^{n_2},$ $n_1+n_2=8;$ $(q^2+1)^{n_1}(q+1)^{n_2}(q-1)^{n_3},$
$1\leqslant
n_1\leqslant4,$ $2n_1+n_2+n_3=8,$ and $|T|\neq(q^2+1)^3(q\pm1)^2,$
$|T|\neq(q^2+1)(q\pm1)^6;$
$(q^3+1)^{n_1}(q^3-1)^{n_2}(q^2+1)^{n_3}(q+1)^{n_4}(q-1)^{n_5},$ $1\leqslant
n_1+n_2\leqslant2,$
$3n_1+3n_2+2n_3+n_4+n_5=8,$ and $|T|\neq(q^3\pm1)^2(q^2+1),$
$|T|\neq(q^3+\epsilon1)(q-\epsilon1)^5,$
$|T|\neq(q^3+\epsilon1)(q^2+1)(q+\epsilon1)^3,$
$|T|\neq(q^3+\epsilon1)(q^2+1)^2(q-\epsilon1);$ $q^8-1;$ $(q^4+1)^2;$
$(q^4+1)(q^2\pm1)(q\pm1)^2;$
$(q^4+1)(q^2-1)^2;$ $(q^4+1)(q^3+\epsilon1)(q-\epsilon1);$
$(q^5+\epsilon1)(q+\epsilon1)^3;$ $(q^5\pm1)(q+\epsilon1)^2(q-\epsilon1);$
$(q^5+\epsilon1)(q^2+1)(q-\epsilon1);$ $(q^5+\epsilon1)(q^3+\epsilon1);$
$(q^6+1)(q^2\pm1);$
$(q^7\pm1)(q\pm1);$ $(q-\epsilon1)\cdot (q^2+\epsilon q+1)^3\cdot(q\pm1);$
$(q^5-\epsilon1)\cdot
  (q^2+\epsilon q +1)\cdot(q+\epsilon1);$ $(q^3\pm1)\cdot
(q^4-q^2+1)\cdot(q\pm1);$
$(q-\epsilon1)\cdot
  (q^6+\epsilon q^3+1)\cdot(q\pm1);$ $(q^3-\epsilon1)\cdot(q^2-\epsilon
  q+1)^2\cdot(q\pm1);$ $q^8-q^4+1;$
  $q^8+q^7-q^5-q^4-q^3+q+1;$
 $q^8-q^6+q^4-q^2+1;$ $(q^4-q^2+1)^2;$ $(q^6+\epsilon q^3+1)(q^2+\epsilon q+1);$
 $q^8-q^7+q^5-q^4+q^3-q+1;$ $(q^4+\epsilon q^3+q^2+\epsilon q+1)^2;$
  $(q^4-q^2+1)(q^2\pm q+1)^2;$
   $(q^2-q+1)^2\cdot(q^2+q+1)^2;$
 $(q^2\pm q+1)^4,$ where $\epsilon=\pm$. Moreover, for every number
given above there exists a torus of corresponding order.
\item[{\em 6.}] Every maximal torus $T$ of $G={^3D_4(q)}$ has one of the
  following orders:
\begin{equation*}
(q^3\pm1)(q\pm1);(q^2\pm q+1)^2;q^4-q^2+1.
\end{equation*}
Moreover, for every number given above there exists a torus of corresponding
order.
\item[{\em 7.}] Every maximal torus $T$ of $G={^2B_2(2^{2n+1})}$ has one of
  the following orders:
\begin{equation*}
q-1; q\pm \sqrt{2q}+1,
\end{equation*}
where $q=2^{2n+1}$. Moreover, for every number given above there exists a torus
of corresponding
order.
\item[{\em 8.}] Every maximal torus $T$ of $G={^2G_2(3^{2n+1})}$ has one of
  the following orders:
\begin{equation*}
q\pm1; q\pm \sqrt{3q} +1,
\end{equation*}
where $q=3^{2n+1}$. Moreover, for every number given above there exists a torus
of corresponding order.
\item[{\em 9.}] Every maximal torus $T$ of $G={^2F_4(2^{2n+1})}$ with $n\ge1$
has one of the following orders: $q^2+\epsilon q\sqrt{2q}+q+\epsilon
\sqrt{2q}+1;$ $q^2-\epsilon
q\sqrt{2q}+\epsilon \sqrt{2q}-1;$ $q^2-q+1;$ $(q\pm \sqrt{2q}+1)^2;$ $(q-1)(q\pm
\sqrt{2q}+1);$
$(q\pm1)^2;q^2\pm1;$ where $q=2^{2n+1}$ and~$\epsilon=\pm$. Moreover, for every
number given above
there exists a torus of corresponding order.
\end{itemize}
\end{lem}

If $q$ is a natural number, $r$ is an odd prime and $(r,q)=1$,
then by $e(r,q)$ we denote the minimal natural number $n$ with
$q^n\equiv1(\mod r)$. If $q$ is odd, let $e(2,q)=1$ if
$q\equiv1(\mod4)$ and $e(2,q)=2$ if $q\equiv-1{(\mod4)}$.

The main technical tools in Section \ref{sporadicalt} is the following
statement.

\begin{lem}\label{Zsigmondy Theorem}
{\em (Corollary to Zsigmondy's theorem
\cite{zs})} Let $q$ be a natural number greater than $1$.
For every natural number $m$ there exists a prime $r$ with $e(r,q)=m$ but for
the cases $q=2$ and
$m=1$, $q=3$ and $m=1$, and $q=2$ and $m=6$.
\end{lem}

The prime $r$ with $e(r,q)=n$ is said to be a {\em primitive prime
divisor} of $q^n-1$. By Zsigmondy theorem it exists excepting the
cases indicated above. If $q$ is fixed, we denote by $r_n$ some
primitive prime divisor of $q^n-1$ (obviously, $q^n-1$ can have more
than one such divisor). Note that according to our definition every
prime divisor of $q-1$ is a primitive prime divisor of $q-1$ with
sole exception: $2$ is not a primitive prime divisor of $q-1$ if
$e(2,q)=2$. In the last case $2$ is a primitive prime divisor
of~$q^2-1$.

In view \cite[Theorems~9.4.10 and~14.3.1]{Ca1} the order of any finite
simple group of Lie type $G$ of rank $n$ over the field $GF(q)$ of
characteristic $p$ is given by $$\vert
G\vert=\frac{1}{d}q^N(q^{m_1}\pm 1)\cdot\ldots\cdot(q^{m_n}\pm 1).$$
It follows that any prime divisor $r$ of $\vert G\vert$ distinct from
the characteristic $p$ is a primitive divisor of $q^m-1$ for some
natural $m$. Thus, the Zsigmondy Theorem allows us to ``find'' prime
divisors of~$\vert G\vert$. Moreover, if $G$ is neither a Suzuki group
nor a Ree group, Lemmas \ref{toriofclassgrps} and
\ref{toriofexcptgrps} imply that for a fixed $m$ every two primitive
prime divisors of $q^m-1$ are adjacent in~$GK(G)$.

For Suzuki and Ree groups we use following

\begin{lem}\label{SuzReeDivisors}
Let $n$ be a natural number.

\noindent {\em 1.} Let
 $m_1(B,n)=2^{2n+1}-1$,

$m_2(B,n)=2^{2n+1}-2^{n+1}+1$,

$m_3(B,n)=2^{2n+1}+2^{n+1}+1$.

Then $(m_i(B,n),m_j(B,n))=1$ if~$i\not=j$.
%\item[{\em 2.}]

\noindent {\em 2.} Let $m_1(G,n)=3^{2n+1}-1$,

$m_2(G,n)=3^{2n+1}+1$,

$m_3(G,n)=3^{2n+1}-3^{n+1}+1$,

$m_4(G,n)=3^{2n+1}+3^{n+1}+1$.

Then $(m_1(G,n), m_2(G,n))=2$ and $(m_i(G,n),m_j(G,n))=1$ otherwise.
%\item[{\em 3.}]

\noindent {\em 3.} Let $m_1(F,n)=2^{2n+1}-1$,

$m_2(F,n)=2^{2n+1}+1$,

$m_3(F,n)=2^{4n+2}+1$,

$m_4(F,n)=2^{4n+2}-2^{2n+1}+1$,

$m_5(F,n)=2^{4n+2}-2^{3n+2}+2^{2n+1}-2^{n+1}+1$,

$m_6(F,n)=2^{4n+2}+2^{3n+2}+2^{2n+1}+2^{n+1}+1$.

Then $(m_2(F,n),m_4(F,n))=3$ and $(m_i(F,n),m_j(F,n))=1$ otherwise.
\end{lem}

\begin{proof}
It easy to check by the direct computation.
\end{proof}

By Lemma~\ref{toriofexcptgrps} every distinct from the characteristic prime
divisor $s$ of order of the
Suzuki group ${}^2B_2(2^{2n+1})$ divides one of the numbers $m_i(B,n)$
defined in Lemma~\ref{SuzReeDivisors}. The same is true for Ree groups
${}^2G_2(3^{2n+1})$, ${}^2F_2(2^{2n+1})$ and all prime divisors of the
numbers $m_i(G,n)$, $m_i(F,n)$ respectively.  Thus,
Lemma~\ref{SuzReeDivisors} allow us to find prime divisors of orders
of Suzuki and Ree groups. Moreover, Lemma~\ref{toriofexcptgrps}
implies that for a fixed $k$ every two prime divisors of $m_k(B,n)$
are adjacent in $GK({}^2B_2(2^{2n+1}))$. The same is also true for Ree
groups and all prime divisors of $m_k(G,n)$ and~$m_k(F,n)$.

\section{Adjacent odd primes}\label{odds}
In this section we consider whether two odd primes distinct from the
characteristic are adjacent in the Gruenberg~--- Kegel graph of a
finite group of Lie type.

\begin{prop}\label{adjanspl}
Let $G=A_{n-1}(q)$ be a finite simple group of Lie type over a
field of characteristic $p$. Let $r,s$ be odd primes and $r,s\in\pi(G)\setminus
\{p\}$. Denote $k=e(r,q)$, $l=e(s,q)$ and suppose that $2\le k\le
l$. Then $r$ and $s$ are non-adjacent if and only if~$k+l> n$ and $k$ does not divide~$l$.
\end{prop}

\begin{proof}
Note first that, for any odd prime $c\ne p$:
\begin{equation}\label{dividers}
c\textit{ divides }
q^x-1\textit{ if and only if }e(c,q)\textit{ divides }x.
\end{equation}
Indeed, by definition, $c$ divides $q^{e(c,q)}-1$ and does
not divide $q^y-1$ for all $y<e(c,q)$, i.~e., $e(c,q)$ is the order of
$q$ in the multiplicative group $GF(c)^\ast$ of the finite field
$GF(c)$.  So if $c$ divides $q^z-1$, then $q^z=1$ in $GF(c)^\ast$,
hence $e(c,q)$ divides $z$. Now assume that $e(c,q)$ divides $z$. Then
$q^z-1=(q^{e(c,q)}-1)f(q)$ for some $f(t)\in \mathbb{Z}[t]$, hence $c$
divides~$q^z-1$.

Assume that $k+l\le n$. Consider a maximal torus $T$
of $G$ of order
$$\frac{1}{(n,q-1)(q-1)} (q^k-1)\cdot(q^l-1)\cdot(q-1)^{n-k-l}.$$ The torus
$T$ is an abelian subgroup of $G$ and $r,s\in\pi(T)$. It follows
that $T$ contains an element of order $rs$, hence $r,s$ are
adjacent. If $k$ divides $l$ then both $r$ and $s$ divide $q^l-1$,
therefore a maximal torus of order
$\frac{1}{(n,q-1)}(q^l-1)(q-1)^{n-l-1}$ contains an element of
order~$rs$.

Assume now that $k+l>n$, $k$ does not divide $l$, and assume that
$g\in G$ is an element of order $rs$. Then $(\vert g\vert, p)=1$
hence $g$ is semisimple. Therefore there exists a maximal torus $T$
such that $g\in T$. By Lemma \ref{toriofclassgrps} the order of $T$
is equal to
$$\frac{1}{(n,q-1)(q-1)} (q^{n_1}-1)\cdot(q^{n_2}-1)\cdot\ldots\cdot
(q^{n_x}-1)$$ for appropriate partition $n_1+n_2+\ldots+n_x=n$ of
$n$. Since $r,s$ are prime, there exist $n_i,n_j$ such that $r$
divides $q^{n_i}-1$ and $s$ divides $q^{n_j}-1$. In view of
\eqref{dividers}, it follows that $n_i=a\cdot k$, $n_j=b\cdot l$ for
some $a,b\ge1$. Moreover, since $k+l>n$ and $k\le l$, we have that
$b=1$. Indeed, otherwise $n_j\ge l+l\ge k+l>n$, a contradiction with
$n_1+\ldots+n_x=n$. Since $k$ does not divide $l$ it follows that
$n_i\not=n_j$. Hence $n_1+n_2+\ldots+n_x\ge n_i+n_j=a\cdot k+ l>n$; a
contradiction.
\end{proof}

\begin{prop}\label{adjantwst}
Let $G={^2A_{n-1}(q)}$ be a finite simple group of Lie type over a
field of characteristic $p$. Define
a function $$\nu(m)=\left\{
\begin{array}{rl}
m &\text{, if }m\equiv 0(\mod 4),\\
\frac{m}{2}& \text{, if }m\equiv 2(\mod 4),\\
2m&\text{, if }m\equiv1(\mod 2).\\
\end{array}\right.$$ Let $r,s$ be odd primes and $r,s\in\pi(G)\setminus
\{p\}$. Denote $k=e(r,q)$, $l=e(s,q)$ and suppose that $2\le \nu(k)\le
\nu(l)$. Then $r$ and $s$ are non-adjacent if and only
if~$\nu(k)+\nu(l)> n$ and $\nu(k)$ does not divide~$\nu(l)$.
\end{prop}

\begin{proof}
Note first that, for any odd prime $c\ne p$:
\begin{equation}\label{dividersantw}
c\textit{ divides
}q^x-(-1)^x\textit{ if and only if }\nu(e(c,q))\textit{ divides
}x.
\end{equation}
Assume first that $e(c,q)$ is odd. If $c$ divides $q^z-(-1)^z$ then
$c$ divides $q^{2z}-1$, hence, by \eqref{dividers}, $e(c,q)$ divides
$2z$. Since $e(c,q)$ is odd it follows that $e(c,q)$ divides $z$ and
therefore (again by \eqref{dividers}) $c$ divides $q^z-1$. Since $c$
is odd we have that $q^z+1$ is not divisible by $c$, so
$q^z-(-1)^z=q^z-1$, i.~e., $z$ is even. But $e(c,q)$ is odd hence
$2e(c,q)=\nu(e(c,q))$ divides $z$. Now assume that $\nu(e(c,q))$
divides $z$. Then $z$ is even, hence $q^z-(-1)^z=q^z-1$ and $c$
divides $q^z-(-1)^z$ by \eqref{dividers}. So \eqref{dividersantw} is
true in this case.

Assume that $e(c,q)\equiv 2(\mod 4)$. If $c$
divides $q^z-(-1)^z$ then $c$ divides $q^{2z}-1$, hence $e(c,q)$
divides $2z$. But $\nu(e(c,q))=\frac{e(c,q)}{2}$, therefore
$\nu(e(c,q))$ divides $z$. If $\nu(e(c,q))$ divides $z$, and $z$ is odd,
then $q^z-(-1)^z=q^z+1$. We have that $e(c,q)$ divides $2z$, hence $c$
divides $q^{2z}-1$. Now $q^{2z}-1=(q^{z}-1)\cdot(q^z+1)$. Since $z$ is
odd, $e(c,q)$ does not divide $z$, therefore, by \eqref{dividers}, $c$
does not divide $q^z-1$, hence $c$ divides $q^z+1=q^z-(-1)^z$. If $\nu(e(c,q))$
divides $z$, and $z$, is even then $2\nu(e(c,q))=e(c,q)$ divides $z$,
hence, by \eqref{dividers}, $c$ divides $q^z-(-1)^z=q^z-1$. Thus
\eqref{dividersantw} is true in this case as well.

At the end assume that
$e(c,q)\equiv 0(\mod 4)$. If $c$ divides $q^z-(-1)^z$ then, as above,
$c$ divides $q^{2z}-1$, hence $e(c,q)$ divides $2z$ and
$\frac{e(c,q)}{2}$ divides $z$. But $\frac{e(c,q)}{2}$ is even, hence,
$z$ is even and $q^z-(-1)^z=q^z-1$. It follows that
$e(c,q)=\nu(e(c,q))$ divides~$z$. If $\nu(e(c,q))=e(c,q)$ divides $z$
then $z$ is even and, by \eqref{dividers}, $c$
divides~$q^z-(-1)^z=q^z-1$.

Assume that $\nu(k)+\nu(l)\le n$. Consider a maximal torus $T$
of $G$ of order
$$\frac{1}{(n,q+1)(q+1)}
(q^{\nu(k)}-(-1)^{\nu(k)})\cdot (q^{\nu(l)}-(-1)^{\nu(l)})\cdot
(q+1)^{n-\nu(k)-\nu(l)}.$$ The torus $T$ is an abelian subgroup of $G$
and $r,s\in\pi(T)$. It follows that $T$ contains an element of order
$rs$, hence $r,s$ are adjacent. If $\nu(k)$ divides $\nu(l)$ then
both $r$ and $s$ divide $q^{\nu(l)}-(-1)^{\nu(l)}$, therefore a
maximal torus of order
$\frac{1}{(n,q+1)}(q^{\nu(l)}-(-1)^{\nu(l)})(q+1)^{n-\nu(l)-1}$
contains an element of order~$rs$.

Assume now that $\nu(k)+\nu(l)>n$, $\nu(k)$ does not divide
$\nu(l)$, and assume that $g\in {^2A_{n-1}(q)}$ is an element of
order $rs$. Then $(\vert g\vert, p)=1$ hence $g$ is semisimple.
Therefore there exists a maximal torus $T$ such that $g\in T$. Using
\eqref{dividersantw} and Lemma \ref{toriofclassgrps} we obtain a
contradiction like in the proof of Proposition~\ref{adjanspl}.
\end{proof}

\begin{prop}\label{adjbn}
Let $G$ be one of simple groups of Lie type,  $B_n(q)$ or $C_n(q)$,  over a
field of characteristic~$p$. Define
$$\eta(m)=\left\{
\begin{array}{cc}
m &\text{ if }m\text{ is odd},\\
\frac{m}{2}& \text{ otherwise}.\\
\end{array}\right.$$ Let $r,s$ be odd primes with $r,s\in\pi(G)\setminus\{p\}$.
Put $k=e(r,q)$ and $l=e(s,q)$, and suppose that $1\le
\eta(k)\le \eta(l)$. Then $r$ and $s$ are non-adjacent if and only if 
$\eta(k)+\eta(l)> n$, and  $k$, $l$ satisfy
to~\eqref{strange}{\em:}
\end{prop}

\begin{equation}\label{strange}
\dfrac{l}{k} \text{ {\em is not an odd natural number}}
\end{equation}

Note that \eqref{strange} is  true in the following cases:

\begin{enumerate}
\item Both $k$ and $l$ are even and either $\eta(k)$ does not divide
  $\eta(l)$ or $\frac{\eta(l)}{\eta(k)}$ is
  even. In this case $q^{\eta(k)}+(-1)^k=q^{\eta(k)}+1$ does not
divide~$q^{\eta(l)}+(-1)^l=q^{\eta(l)}+1$.
\item Both $k$ and $l$ are odd and $\eta(k)$ does not divide
  $\eta(l)$. In this case $q^{\eta(k)}+(-1)^k=q^k-1$ does not
  divide~$q^{\eta(l)}+(-1)^l=q^l-1$.
\item $k$ is odd and $l$ is even. In this case  $q^{\eta(k)}+(-1)^k=q^k-1$ does not
  divide~$q^{\eta(l)}+(-1)^l=q^{\eta(l)}+1$.
\item $k$ is even, $l$ is odd and either $\eta(k)$ does not divide
  $\eta(l)$ or $\frac{\eta(l)}{\eta(k)}$ is either even, or is equal
  to 1. In this case $q^{\eta(k)}+(-1)^k=q^{\eta(k)}+1$ does not
  divide~$q^{\eta(l)}+(-1)^l=q^l-1$.
\end{enumerate}

\noindent This means that $\eta(k)$, $\eta(l)$ satisfy
\eqref{strange} if and only if $q^{\eta(k)}+(-1)^k$ does not
divide~${q^{\eta(l)}+(-1)^l}$.

\begin{proof}
First we prove that, for any odd prime $c\ne p$:
\begin{equation}\label{dividersbn}
\textit{if }c\textit{ divides
}q^x\pm1\textit{ then }\eta(e(c,q))\textit{
divides }x.
\end{equation}
\noindent Assume first that $e(c,q)$ is odd and $c$ divides
$q^z\pm1$. By \eqref{dividers} it follows that
$e(c,q)$ divides $2z$. But $e(c,q)$ is odd hence
$e(c,q)=\eta(e(c,q))$ divides $z$. Assume now that $e(c,q)$ is even and $c$ divides
$q^z\pm1$. It follows that $c$ divides $q^{2z}-1$, hence $e(c,q)$
divides $2z$ and $\frac{e(c,q)}{2}=\eta(e(c,q))$ divides~$z$.

Now if $\eta(k)+\eta(l)\le n$ we may consider a maximal torus $T$ of
order $\frac{1}{(2,q-1)}(q^{\eta(k)}+(-1)^k)\cdot
(q^{\eta(l)}+(-1)^l)\cdot (q-1)^{n-\eta(k)-\eta(l)}$. We have that
$T$ is an abelian group and $r,s\in\pi(T)$. Hence $T$ contains an
element of order~$rs$. If $\eta(k)$, $\eta(l)$ does not satisfy
\eqref{strange}, then $q^{\eta(k)}+(-1)^k$ divides
$q^{\eta(l)}+(-1)^l$, therefore both $r$ and $s$ divide
$q^{\eta(l)}+(-1)^l$. Thus  a maximal torus of order
$\frac{1}{(2,q-1)}(q^{\eta(l)}+(-1)^l)(q-1)^{n-\eta(l)}$ contains an
element of order~$rs$.

Assume that $\eta(k)+\eta(l)>n$, $\eta(k)$ and $\eta(l)$ satisfy
\eqref{strange}, and assume that there exists an element $g\in G$ of
order $rs$. Since $(\vert g\vert,p)=1$, it follows that $g$ is
semisimple. So there exists a maximal torus $T$ containing $g$. In
view of Lemma \ref{toriofclassgrps}, the order $\vert T\vert$ is
equal to
$$\frac{1}{(2,q-1)}
(q^{n_1}\pm1)\cdot(q^{n_2}\pm1)\cdot\ldots\cdot(q^{n_x}\pm1)$$ for
appropriate partition $n_1+n_2+\ldots+n_x=n$ of
$n$.  Using
\eqref{dividersbn} we obtain a contradiction like in
Proposition~\ref{adjanspl}.
\end{proof}

\begin{prop}\label{adjdn}
et $G=D_n^{\varepsilon}(q)$ be a finite simple group of Lie type over a field of
characteristic  $p$, and let the function
$\eta(m)$ be defined as in Proposition~{\em\ref{adjbn}}. Suppose $r,s$ are odd
primes and $r,s\in\pi(D_n^\varepsilon(q))\setminus\{p\}$.
Put $k=e(r,q)$, $l=e(s,q)$, and  $1\le\eta(k)\le\eta(l)$. Then $r$ and $s$ are
non-adjacent if and only if $2\cdot\eta(k)+2\cdot\eta(l)>
2n-(1-\varepsilon(-1)^{k+l})$,
$k$ and $l$ satisfy~\eqref{strange}, and, if  $\varepsilon=+$, then the chain of
equalities{\em:}
\begin{equation}\label{strange2}
n=l=2\eta(l)=2\eta(k)=2k
\end{equation}
is not true.
\end{prop}

\begin{proof}
Using \eqref{dividersbn} and Lemma \ref{toriofclassgrps}   we prove
the proposition as above.
\end{proof}

\begin{prop}\label{adjexcept}
Let $G$ be a finite simple exceptional group of Lie type over a field of
characteristic~$p$, suppose that
$r,s$ are odd primes, and assume that  $r,s\in\pi(G)\setminus\{p\}$, $k=e(r,q)$,
$l=e(s,q)$, and
$1\le k\le l$. Then  $r$ and $s$ are non-adjacent if and only if  $k\not=l$ and
one of the following holds{\em:}
\begin{itemize}
\item[{\em 1.}] $G=G_2(q)$  and either $r\not=3$ and $l\in\{3,6\}$ or $r=3$
and~${l=9-3k}$.
\item[{\em 2.}] $G=F_4(q)$ and either $l\in\{8,12\}$, or $l=6$
and~$k\in\{3,4\}$, or $l=4$ and~${k=3}$.
\item[{\em 3.}] $G=E_6(q)$ and either $l=4$ and $k=3$, or $l=5$ and $k\ge3$, or
$l=6$
and $k=5$, or $l=8$, $k\ge3$, or $l=8$, $r=3$, and $(q-1)_3=3$, or
$l=9$, or $l=12$ and $k\not=3$.
\item[{\em 4.}] $G={^2E_6(q)}$ and either $l=6$ and $k=4$, or $l=8$, $k\ge3$, or
$l=8$, $r=3$,and
$(q+1)_3=3$, or  $l=10$ and $k\ge 3$, or $l=12$ and $k\not=6$, or~$l=18$.
\item[{\em 5.}] $G=E_7(q)$ and either $l=5$ and $k=4$, or  $l=6$ and $k=5$, or
$l\in\{14,18\}$ and
$k\not=2$, or  $l\in\{7,9\}$ and $k\ge2$, or $l=8$ and $k\ge3,k\not=4$,  or
$l=10$ and $k\ge3,
  k\not=6$, or $l=12$ and~$k\ge 4,k\not=6$.
\item[{\em 6.}] $G=E_8(q)$ and either $l=6$ and $k=5$, or $l\in\{7,14\}$ and
$k\ge3$, or
$l=9$ and $k\ge 4$, or $l\in\{8,12\}$ and $k\ge 5, k\not=6$, or $l=10$ and
$k\ge3, k\not=4,6$, or $l=18$ and $k\not=1,2,6$, or
$l=20$ and $r\cdot k\not=20$, or~$l\in\{15,24,30\}$.
\item[{\em 7.}] $G={^3D_4(q)}$ and either $l=6$ and $k=3$,
or~$l=12$.
\end{itemize}
\end{prop}

\begin{proof}
As for classical groups of Lie type, prime divisors $r,s\in\pi(G)$
satisfying the conditions of the proposition are adjacent if and
only if $rs$ divides the order of some maximal torus of $G$. Thus,
using Lemma~\ref{toriofexcptgrps} instead of
Lemma~\ref{toriofclassgrps} we prove the proposition as above.
\end{proof}

\begin{prop}\label{adjsuzree}
Let $G$ be a finite simple Suzuki or Ree group over a field of
characteristic $p$, let $r,s$  be odd primes $r,s\in\pi(G)\setminus\{p\}$. Then
$r,s$ are non-adjacent if and only if one of the
following holds:
\begin{itemize}
\item[{\em 1.}] $G={^2B_2(2^{2n+1})}$, $r$ divides $m_k(B,n)$, $s$ divides $m_l(B,n)$ and~$k\not=l$.
\item[{\em 2.}] $G={^2G_2(3^{2n+1})}$, $r$ divides $m_k(G,n)$, $s$ divides $m_l(G,n)$ and~$k\not=l$.
\item[{\em 3.}] $G={^2F_4(2^{2n+1})}$, $r$ divides $m_k(F,n)$, $s$
  divides $m_l(F,n)$, $k\not=l$, and if $\{k,l\}=\{1,3\}$,
then~${r\not=3\not=s}$.
\end{itemize}
Numbers $m_i(B,n)$, $m_i(G,n)$, and $m_i(F,n)$ are defined in
Lemma~{\em\ref{SuzReeDivisors}}.
\end{prop}

\begin{proof}
We use Lemma \ref{toriofexcptgrps}, Lemma \ref{SuzReeDivisors} and
arguments as in the previous propositions of the section.
\end{proof}

\section{Adjacency with the characteristic}\label{char}

In this section we consider whether a prime $r$ and the
characteristic $p$ of the base field of a finite group of Lie type
are adjacent.

\begin{prop}\label{adjcharclass}
Let $G=O^{p'}(\overline{G}_\sigma)$ be a finite simple classical
group of Lie type defined over a field of characteristic $p$. Let
$r\in\pi(G)$ and $r\not=p$. Then $r$ and $p$ are non-adjacent if and only
if one of the following holds:

\begin{itemize}
\item[{\em 1.}] $G=A_{n-1}(q)$, $r$ is odd, and $e(r,q)> n-2$.
\item[{\em 2.}] $G={^2A_{n-1}(q)}$, $r$ is odd, and $\nu(e(r,q))> n-2$. The
function $\nu(m)$ is defined in Proposition~{\em\ref{adjantwst}}.
\item[{\em 3.}] $G=C_n(q)$,  $\eta(e(r,q))> n-1$. The function
$\eta(m)$ is defined in Proposition~{\em\ref{adjbn}}.
\item[{\em 4.}] $G=B_n(q)$,   $\eta(e(r,q))> n-1$. The
function $\eta(m)$ is defined in Proposition~{\em\ref{adjbn}}.
\item[{\em 5.}] $G=D_n^\varepsilon(q)$, $\eta(e(r,q))> n-2$. The function $\eta(m)$ is
defined in Proposition~{\em\ref{adjbn}}.
\item[{\em 6.}] $G=A_1(q)$, $r=2$.
\item[{\em 7.}] $G=A_2^\varepsilon(q)$, $r=3$ and~$(q-\varepsilon1)_3=3$.
\end{itemize}
\end{prop}

\begin{proof}
It is evident that $2$ and $p$ are adjacent in all classical simple groups
except $A_1(q)$. Hence we may assume that $r$ is odd.

First we outline the general idea of the proof. We shall use
\cite{Ca2} as a technical tool in our proof, so we shall keep
notations of \cite{Ca2}. In order to prove that $r$ and $p$ are
adjacent we find a connected reductive subgroup $R$ of maximal rank
in $G$ such that $R=T(G_1\ast G_2)$, where $T$ is a maximal torus,
both $G_1$ and $G_2$ are non-trivial groups of Lie type such that
$r\in\pi(G_1)$. Then $G_1$ contains an element of order $r$ and it
centralizes $G_2$. Since $G_2$ is non-trivial, it contains an
element of order $p$, so $G\geq R$ contains an element of
order~$rp$.

In order to prove that $r$ and $p$ are non-adjacent we shall
consider arbitrary element $g$ of order $r$ and its connected
centralizer $G\cap C_{\overline{G}}(g)^0$. Recall that
$C_{\overline{G}}(g)^0=\overline{S}\ast \overline{L}$, where
$\overline{S}=Z(C_{\overline{G}}(g)^0)$ is a central torus and
$\overline{L}$ is a semisimple part. Clearly $g\in \overline{S}\cap
G\leq \overline{S}_\sigma$, hence $r$ divides $\vert
\overline{S}_\sigma\vert$. This condition would imply that
$\overline{L}_\sigma$ is trivial, hence $C_{\overline{G}}(g)^0$ does
not contain unipotent elements. But every unipotent element of
$C_{\overline{G}}(g)$ is contained in $C_{\overline{G}}(g)^0\leq
N_{\overline{G}}(\overline{T})^0=\overline{T}$. So none unipotent
element of $G$ centralizes~$g$. Now consider all classical groups
case by case.

$A_1(q)$. It is known that if $g$ is an element of order $r\not=p$,
then $(\vert C_{A_1(q)}(g)\vert,p)=1$ (see,
\cite[Proposition~7]{Ca2}). So $r,p$ are non-adjacent for
every~$r\in\pi(A_1(q))\setminus\{p\}$.

$A_2(q)$. By using \cite[Proposition~7]{Ca2} we obtain that only prime
divisors of $\frac{q-1}{(3,q-1)}$ are adjacent to $p$ and the
proposition is true in this case.

$A_{n-1}(q)$, and $n\ge4$. In our case $T$ is a Cartan subgroup,
$G_1=A_{n-3}(q)$ and $G_2=A_1(q)$. The
existence of such a subgroup can be obtained by using
\cite[Proposition~7]{Ca2}. Thus every $r$ with $e(r,q)\le n-2$ divides
$\vert G_1\vert$, hence is adjacent to $p$. Now let $e(r,q)=n-1$ and let
$g$ be an element of order $r$. In view of \eqref{dividers} we obtain
that $q^{n-1}-1$ must divide $\vert
\overline{S}_\sigma\vert\cdot(q-1)$. It follows from \cite[Proposition~7]{Ca2}
that $G\cap C_{\overline{G}}(g)^0$ is a maximal torus of order
$\frac{1}{(n,q-1)}(q^{n-1}-1)$. Hence $\vert
\overline{L}_\sigma\vert=1$ and $r,p$ are non-adjacent. If $e(r,q)=n$
and $g$ is an element of order $r$, then by \eqref{dividers},
$q^n-1$ must divide $\vert
\overline{S}_\sigma\vert\cdot(q-1)$. Again by using \cite[Proposition~7]{Ca2} we
obtain that $G\cap C_{\overline{G}}(g)^0$ is a maximal torus of order
$\frac{1}{(n,q-1)}\frac{q^n-1}{q-1}$. So $\vert
\overline{L}_\sigma\vert=1$ in this case and every $r$ with $e(r,q)>
n-2$ is non-adjacent to~$p$.

${}^2A_2(q)$. By using \cite[Proposition~8]{Ca2} we obtain that only prime
divisors of $\frac{q+1}{(3,q+1)}$ are adjacent to $p$ and the
proposition is true in this case.

${^2A_{n-1}(q)}$ and $n\ge4$. Again $T$ is a Cartan subgroup,
$G_1={^2A_{n-3}(q)}$ and $G_2=A_1(q)$. The existence of such a subgroup can be obtained by using
\cite[Proposition~8]{Ca2}. Thus every $r$ with $\nu(e(r,q))\le n-2$ divides
$\vert G_1\vert$, hence is adjacent to $p$. Now let $\nu(e(r,q))=n-1$
and $g$ is an element of order $r$. Then in view of
\eqref{dividersantw} we obtain that
$q^{n-1}-(-1)^{n-1}$ divides $\vert \overline S\vert\cdot(q+1)$. It
follows from \cite[Proposition~8]{Ca2} that $G\cap
C_{\overline{G}}(g)^0$ is a maximal torus of order
$\frac{1}{(n,q+1)}(q^{n-1}-(-1)^{n-1})$. Hence
$\vert\overline{L}_\sigma\vert=1$ and $r,p$ are non-adjacent. If
$\nu(e(r,q))=n$ and $g$ is an element of order $r$, then by
\eqref{dividersantw}, $q^n-(-1)^n$ divides $\vert
\overline{S}_\sigma\vert\cdot (q+1)$. By using
\cite[Proposition~8]{Ca2} we obtain that $G\cap C_{\overline{G}}(g)^0$
is a maximal torus of order
$\frac{1}{(n,q+1)}\frac{q^n-(-1)^n}{q+1}$. Therefore
$\vert\overline{L}_\sigma\vert=1$ and every $r$ with $\nu(e(r,q))>
n-2$ is non-adjacent to~$p$.

$C_n(q)$. Take for $R$ torus $T$ as a Cartan subgroup,
$G_1=C_{n-1}(q)$, $G_2=A_1(q)$. Such a subgroup $R$ exists in view
of \cite[Propositions~9 and~12]{Ca2}. Again every $r$ with
$\eta(e(r,q))\le n-1$ divides $\vert G_1\vert$, hence is adjacent to
$p$. If $\eta(e(r,q))=n$ and $g$ is an element of order $r$, then
\eqref{dividersbn} implies that either $q^n-1$, or $q^n+1$ divides
$\vert \overline{S}_\sigma\vert$. From \cite[Propositions~9
and~12]{Ca2} we obtain that $G\cap C_{\overline{G}}(g)^0$ is a
maximal torus of order $\frac{1}{(2,q-1)}(q^n\pm1)$ respectively,
hence $\vert\overline{L}_\sigma\vert=1$ and $r,p$ are non-adjacent.

$B_n(q)$. Since $B_2(q)\simeq C_2(q)$ and $B_n(2^t)\simeq C_n(2^t)$,
we may assume that $p$ is odd and $n\ge3$. We can take $T$ to be a
Cartan subgroup, $G_1=B_{n-2}(q)$, and $G_2=D_2(q)$. Such a subgroup
exists in view of \cite[Proposition~11]{Ca2}. Since every prime $r$
with $\eta (e(r,q))\le n-2$ divides the order of $G_1$, we obtain that
$r$ and $p$ are adjacent if $\eta(e(r,q))\le n-2$. If
$\eta(e(r,q))=n-1$, then, again by \cite[Proposition~11]{Ca2}, there
exists a reductive subgroup $R$ such that
$\vert\overline{S}_\sigma\vert=q^{\eta(e(r,q))}+(-1)^{e(r,q)}$ and
$\overline{L}_\sigma=B_1(q)\simeq A_1(q)$. Hence $r,p$ are adjacent in
$G$. Assume now that $\eta(e(r,q))= n$ and that $g$ is an element of
order $r$ of $G$. Then the order $\vert\overline{S}_\sigma\vert$ is
given in \cite[Proposition~11]{Ca2} and is equal to
$\prod_{i}(q^{n_i}\pm1)$, where $\sum_i n_i\le n$. By using
\eqref{dividersbn} we obtain that $\eta(e(r,q))=n$ divides $n_i$ for
some $i$. Hence either $q^n-1$, or $q^n+1$ divides $\vert
S_\sigma\vert$. In view of \cite[Proposition~11]{Ca2} this implies
that~$\vert \overline{L}_\sigma\vert=1$, hence $r,p$ are non-adjacent.

$D_n^\varepsilon(q)$. We can take $T$ to be a Cartan subgroup,
$G_1=D_{n-2}^\varepsilon(q)$, $G_2=A_1(q)$. The existence of such a
subgroup $R$ can be obtained by using \cite[Proposition~10]{Ca2}. So
$r,p$ are adjacent for all $r$ with $\eta(e(r,q))\le n-2$, except one
case: $\eta(e(r,q))=n-2$ and $r$ divides $q^{n-2}+\varepsilon1$. In this
last case we obtain that there exists a reductive subgroup $R$ such
that $\vert\overline{S}_\sigma\vert=q^{n-2}+\varepsilon1$ and
$\overline{L}_\sigma\simeq D_2^\varepsilon(q)$. Hence in this
exceptional case $r,p$ are adjacent also. Now if $\eta(e(r,q))\ge
n-1$ and $g$ is an element of $G$ of order $r$, then as above we
obtain that $q^{\eta(e(r,q))}+(-1)^{e(r,q)}$ divides $\vert
\overline{S}_\sigma\vert$ and \cite[Proposition~10]{Ca2} implies that
$\overline{L}_\sigma$ is trivial. Therefore $r,p$ are non-adjacent in
this case.
\end{proof}

\begin{prop}\label{adjcharexcept}
Let $G$ be a finite simple exceptional group of Lie type over a
field of characteristic $p$. Let $r\in\pi(G)$, $k=e(r,q)$, and
$r\not=p$. Then $r,p$ are non-adjacent if and only if one of the
following holds:

\begin{itemize}
\item[{\em 1.}] $G=G_2(q)$, $k\in\{3,6\}$.
\item[{\em 2.}] $G=F_4(q)$, $k\in\{8,12\}$.
\item[{\em 3.}] $G=E_6(q)$, $k\in\{8,9,12\}$.
\item[{\em 4.}] $G={}^2E_6(q)$, $k\in\{8,12,18\}$.
\item[{\em 5.}] $G=E_7(q)$, $k\in\{7,9,14,18\}$.
\item[{\em 6.}] $G=E_8(q)$, $k\in\{15,20,24,30\}$.
\item[{\em 7.}] $G={^3D_4(q)}$, $k=12$.
\end{itemize}
\end{prop}

\begin{proof}
All statements  are obtained by using information about conjugacy
classes and centralizers of semisimple elements given in \cite{der1}
and~\cite{der2}.
\end{proof}

\begin{prop}\label{adjcharsuzree}
Let $G$ be a finite simple Suzuki or Ree group over a field of
characteristic $p$, let $r\in\pi(G)\setminus\{p\}$. Then $r,p$ are
non-adjacent if and only if one of the following holds:
\begin{itemize}
\item[{\em 1.}] $G={^2B_2(2^{2n+1})}$, $r$ divides $m_k(B,n)$.
\item[{\em 2.}] $G={^2G_2(3^{2n+1})}$, $r$ divides $m_k(G,n)$ and $r\not=2$.
\item[{\em 3.}] $G={^2F_4(2^{2n+1})}$, $r$ divides $m_k(F,n)$,  $r\not=3$, and
$k>2$.
\end{itemize}
Numbers $m_i(B,n)$, $m_i(G,n)$, and $m_i(F,n)$ are defined in
Lemma~{\em\ref{SuzReeDivisors}}.
\end{prop}

\begin{proof}
All statements  are obtained by using information about conjugacy
classes and centralizers of semisimple elements given
in~\cite{der1}.
\end{proof}

\section{Adjacency of 2 and an odd prime $r$}\label{two}

In this section we consider whether $2$ and an odd prime $r$ are
adjacent if both of them are not equal to the characteristic $p$ of
the base field. We start from groups $A_n^\varepsilon(q)$. Recall
that for those groups we did not consider adjacency criterion for
prime divisors of $q-\varepsilon1$ in Section \ref{odds}. It is
natural to consider such a criterion together with a criterion for
2, hence we give it in the following two propositions.

\begin{prop}\label{bigprimesan}
Let $G=A_{n-1}(q)$ be a finite simple group of Lie type. Let $r$ be
a prime divisor of $q-1$ and $s$ be an odd prime distinct from the
characteristic. Denote $k=e(s,q)$. Then $s$ and $r$ are non-adjacent if
and only if one of the following holds:
\begin{itemize}
\item[{\em 1.}] $k=n$, $n_r\le(q-1)_r$, and if $n_r=(q-1)_r$, then $2<(q-1)_r$.
\item[{\em 2.}] $k=n-1$ and~$(q-1)_r\le n_r$.
\end{itemize}
\end{prop}

\begin{proof}
First we prove the following statement

\begin{multline}\label{divofsum}
\left(\frac{q^n-1}{q-1}\right)_r=n_r\text{, if }(q-1)_r\ge n_r\text{
  and }(q-1)_r>2,\\
\left(\frac{q^n-1}{q-1}\right)_2>2\text{, if }(q-1)_2=n_2=2.
\end{multline}

Assume that $r$ is odd. Then $n=r^k\cdot l$, where $(l,r)=1$ and
$q=r^k\cdot m+1$. Now $q^n-1=(q^{r^k}-1)\cdot
(q^{n-r^k}+q^{n-2r^k}+\ldots+q^{r^k}+1)$. The second multiplier is
the sum of $l$ numbers of type $q^i$ and $q^i\equiv 1(\mod r)$ for
all $i$. Since $(l,r)=1$ it follows that the second multiplier is
coprime to $r$. Thus
$\left(\frac{q^n-1}{q-1}\right)_r=\left(\frac{q^{r^k}-1}{q-1}\right)_r$.
Now remember that $q=r^k\cdot m+1$, hence
$$\frac{q^{r^k}-1}{q-1}=(r^k\cdot m)^{r^k-1}+r^k\cdot(r^k\cdot
m)^{r^k-2}+\ldots+\frac{1}{2}r^k(r^k-1)(r^k\cdot m)+r^k.$$ Since
$r^{k+1}$ divides all summands, except the last, and $r^{k+1}$ does
not divide $r^k$ we obtain that $r^k$ divides
$\left(\frac{q^{r^k}-1}{q-1}\right)_r$, but $r^{k+1}$ does not divide
$\left(\frac{q^{r^k}-1}{q-1}\right)_r$. Therefore
$\left(\frac{q^{r^k}-1}{q-1}\right)_r=n_r$ in this case.

Assume now that $r=2$, $n=2^k\cdot l$, where $l$ is odd, and
$q=2^k\cdot m+1$. Then
$$\left(\frac{q^n-1}{q-1}\right)_2=\left(\frac{q^l-1}{q-1}\right)_2\cdot
(q^l+1)\cdot (q^{2l}+1)\cdot \ldots\cdot (q^{2^{k-1}l}+1).$$ Since
$l$ is odd we have that $\left(\frac{q^l-1}{q-1}\right)_2=1$. If
$(q-1)_2>2$ we have that $(q^i+1)_2=2$ for all $i$. Therefore
$$\left(\frac{q^n-1}{q-1}\right)_2= (q^l+1)_2\cdot (q^{2l}+1)_2\cdot
\ldots\cdot (q^{2^{k-1}l}+1)_2 =2^k=n_2.$$ The second case
$(q-1)_2=n_2=2$ is evident.

Assume that $k\le n-2$. By Lemma \ref{toriofclassgrps} it follows
that there exists a maximal torus $T$ of $G$ of order
$\frac{1}{(n,q-1)}(q^k-1)(q-1)^{n-k-1}$. We have that $T$ is an
abelian subgroup of $G$, and both $r,s$ divide $\vert T\vert$.
Therefore $T$ contains an element of order~$rs$.

Assume that $k=n$. By \eqref{dividers} and Lemma
\ref{toriofclassgrps} every element of order $s$ is contained in a
maximal torus $T$ of order $\frac{1}{(n,q-1)}\frac{q^n-1}{q-1}$. In
view of \eqref{divofsum} $r$ does not divide $\vert T\vert$ if and
only if condition 1 of the proposition holds.

Assume now that $k=n-1$. By \eqref{dividers} and Lemma
\ref{toriofclassgrps} every element of order $s$ is contained in a
maximal torus $T$ of order $\frac{1}{(n,q-1)}(q^{n-1}-1)$. We have
that $(q^{n-1}-1)=(q-1)(q^{n-2}+q^{n-3}+\ldots+q+1)$ and
$(q^{n-2}+q^{n-3}+\ldots+q+1)_r=1$. Therefore $r$ does not divide
$\vert T\vert$ if and only if $(q-1)_r\le n_r$ and we obtain
condition 2 of the proposition.
\end{proof}

\begin{prop}\label{bigprimesantw}
Let $G={^2A_{n-1}(q)}$ be a finite simple group of Lie type. Let $r$
be a prime divisor of $q+1$ and $s$ be an odd prime distinct from
the characteristic. Denote $k=e(s,q)$. Then $s$ and $r$ are non-adjacent
if and only if one of the following holds:
\begin{itemize}
\item[{\em 1.}] $\nu(k)=n$, $n_r\le(q+1)_r$, and if $n_r=(q+1)_r$, then $2<(q+1)_r$.
\item[{\em 2.}] $\nu(k)=n-1$ and~$(q+1)_r\le n_r$.
\end{itemize}
The function $\nu(m)$ is defined in Proposition~{\em\ref{adjantwst}}.
\end{prop}

\begin{proof}
Like in the previous proposition first we show that

\begin{multline}\label{divofsumtw}
\left(\frac{q^n-(-1)^n}{q+1}\right)_r=n_r\text{, if }(q+1)_r\ge n_r\text{
  and }(q+1)_r>2,\\
\left(\frac{q^n-(-1)^n}{q+1}\right)_2>2\text{, if }(q+1)_2=n_2=2.
\end{multline}

Assume that $r$ is odd. Then $n=r^k\cdot l$, where $(l,r)=1$ and
$q=r^k\cdot m-1$. Now $q^n-(-1)^n=(q^{r^k}+1)\cdot
(q^{n-r^k}-q^{n-2r^k}+\ldots+ (-1)^{l-1}q^{r^k}+(-1)^l)$. The second
multiplier is the sum of $l$ numbers of type $(-1)^t q^{n-(t+1)r^k}$
and $q^{n-(t+1)r^k}\equiv (-1)^{n+t-1}(\mod r)$ for all $t$. Hence
the second multiplier is equivalent $(-1)^{n-1}l$ modulo $r$. Since
$(l,r)=1$ it follows that the second multiplier is coprime to $r$.
Thus
$\left(\frac{q^n-(-1)^n}{q+1}\right)_r=\left(\frac{q^{r^k}+1}{q+1}\right)_r$.
Now remember that $q=r^k\cdot m-1$, hence
$$\frac{q^{r^k}+1}{q+1}=(r^k\cdot m)^{r^k-1}-r^k\cdot(r^k\cdot
m)^{r^k-2}+\ldots+(-1)^{k-2}\frac{1}{2}r^k(r^k-1)(r^k\cdot m)+(-1)^{k-1}r^k.$$ Since
$r^{k+1}$ divides all summands, except the last, and $r^{k+1}$ does
not divide $r^k$ we obtain that $r^k$ divides
$\left(\frac{q^{r^k}+1}{q+1}\right)_r$, but $r^{k+1}$ does not divide
$\left(\frac{q^{r^k}+1}{q+1}\right)_r$. Therefore
$\left(\frac{q^{r^k}+1}{q+1}\right)_r=n_r$ in this case.

Assume now that $r=2$, $n=2^k\cdot l$, where $l$ is odd, and
$q=2^k\cdot m-1$. Then
$$\left(\frac{q^n-(-1)^n}{q+1}\right)_2=\left(\frac{q^l+1}{q+1}\right)_2\cdot
(q^l-1)\cdot (q^{2l}+1)\cdot \ldots\cdot (q^{2^{k-1}l}+1).$$ Since
$l$ is odd we have that $\left(\frac{q^l+1}{q+1}\right)_2=1$. If
$(q+1)_2>2$ we have that $(q^{2i}+1)_2=2$ for all $i\ge1$ and
$(q^l-1)_2=2$. Therefore
$$\left(\frac{q^n-(-1)^n}{q+1}\right)_2= (q^l-1)_2\cdot (q^{2l}+1)_2\cdot
\ldots\cdot (q^{2^{k-1}l}+1)_2 =2^k=n_2.$$ The second case
$(q+1)_2=n_2=2$ is evident.

Assume that $\nu(k)\le n-2$. By Lemma \ref{toriofclassgrps} it
follows that there exists a maximal torus $T$ of $G$ of order
$\frac{1}{(n,q+1)}(q^{\eta(k)}-(-1)^{\eta(k)})(q+1)^{n-k-1}$. We
have that $T$ is an abelian subgroup of $G$ and both $r,s$ divide
$\vert T\vert$. Therefore $T$ contains an element of order~$rs$.

Assume that $\nu(k)=n$. By \eqref{dividersantw} and Lemma
\ref{toriofclassgrps} every element of order $s$ is contained in a
maximal torus $T$ of order
$\frac{1}{(n,q+1)}\frac{q^n-(-1)^n}{q+1}$. In view of
\eqref{divofsumtw} $r$ does not divide $\vert T\vert$ if and only if
condition 1 of the proposition holds.

Assume now that $\nu(k)=n-1$. By \eqref{dividersantw} and
Lemma~\ref{toriofclassgrps} every element of order $s$ is contained
in a maximal torus $T$ of order
$\frac{1}{(n,q+1)}(q^{n-1}-(-1)^{n-1})$. We have that
$(q^{n-1}-(-1)^{n-1})=(q+1)(q^{n-2}-q^{n-3}+\ldots+(-1)^{n-3}q+(-1)^{n-2})$
and $(q^{n-2}-q^{n-3}+\ldots+(-1)^{n-3}q+(-1)^{n-2})_r=1$. Therefore
$r$ does not divide $\vert T\vert$ if and only if $(q+1)_r\le n_r$
and we obtain condition 2 of the proposition.
\end{proof}

\begin{prop}\label{adj2bn}
Let $G=B_n(q)$ or $G=C_n(q)$ be a finite simple group of Lie type over
a field of odd characteristic $p$. Let $r$ be an odd prime divisor of
$\vert G\vert$, $r\not=p$, and $k=e(r,q)$. Then $2,r$ are non-adjacent if and
only
if $\eta(k)=n$ and one of the following holds:
\begin{itemize}
\item[{\em 1.}] $n$ is odd and~$k=(3-e(2,q))n$.
\item[{\em 2.}] $n$ is even and~$k=2n$.
\end{itemize}
The function $\eta(m)$ is defined in Proposition~{\em\ref{adjbn}}.
\end{prop}

\begin{proof}
If $\eta(k)\le n-1$, then by Lemma \ref{toriofclassgrps} there
exists a maximal torus $T$ of order
$\frac{1}{2}(q^{\eta(k)}+(-1)^k)(q-1)^{n-\eta(k)}$. We have that $T$
is an abelian subgroup of $G$ and $2,r$ both divide $\vert T\vert$.
Hence $T$ contains an element of order $2r$ and $2,r$ are adjacent.

Thus we may assume that $\eta(k)=n$. In view of \eqref{dividersbn}
and Lemma \ref{toriofclassgrps} we have that every element $g$ of
order $r$ is contained in a maximal torus $T$ of order
$\frac{1}{2}(q^n+(-1)^k)$. Therefore $2,r$ are non-adjacent if and
only if 2 does not divide $\vert T\vert$ and the proposition
follows.
\end{proof}

\begin{prop}\label{adj2dn}
Let $G=D_n^\varepsilon(q)$ be a finite simple group of Lie type over a field of odd
characteristic $p$. Let $r$ be an odd prime divisor of $\vert G\vert$,
$r\not=p$, and $k=e(r,q)$. Then $2$ and $r$ are non-adjacent if and only if one
of the
following holds:
\begin{itemize}
\item[{\em 1.}] $\eta(k)=n$
  and~$(4,q^n-\varepsilon1)=(q^n-\varepsilon1)_2$.
\item[{\em 2.}] $\eta(k)=k=n-1$, $n$ is even, $\varepsilon=+$,
and~$e(2,q)=2$.
\item[{\em 3.}] $\eta(k)=\frac{k}{2}=n-1$, $\varepsilon=+$
  and~$e(2,q)=1$.
\item[{\em 4.}]   $\eta(k)=\frac{k}{2}=n-1$, $n$ is odd,
  $\varepsilon=-$, and~$e(2,q)=2$.
\end{itemize}
The function $\eta(m)$ is defined in Proposition~{\em\ref{adjbn}}.
\end{prop}

\begin{proof}
First note that if $\eta(k)\le n-2$, then by Lemma
\ref{toriofclassgrps} there exists a maximal torus $T$ of order
$\frac{1}{(4,q^n-\varepsilon1)}(q^{\eta(k)}+(-1)^k)(q+(-\varepsilon1)^k)
(q-1)^{n-\eta(k)-1}$. We have that $2,r$ both divide $\vert T\vert$,
hence $2,r$ are adjacent.

Now consider the case $G={^2D_n(q)}$ (i.~e. $\varepsilon=-$). If
$\eta(k)=n-1$ then, by \eqref{dividersbn} and Lemma
\ref{toriofclassgrps}, we have that every element of order $r$ is
contained in a maximal torus $T$ of order
$\frac{1}{(4,q^n+1)}(q^{n-1}+(-1)^k)(q-(-1)^k)$. Then 2 does not
divide $\vert T\vert$ if and only if condition~4 of the proposition
holds. If $\eta(k)=n$, then, by \eqref{dividersbn} and Lemma
\ref{toriofclassgrps}, every element of order $r$ is contained in a
maximal torus $T$ of order $\frac{1}{(4,q^n+1)}(q^n+1)$ (in
particular, $k=2\eta(k)$). Then 2 does not divide $\vert T\vert$ if
and only if $(4,q^n+1)=(q^n+1)_2$ and we obtain the proposition for twisted
groups.

Assume now that $G=D_n(q)$. If $\eta(k)=n-1$, then, by
\eqref{dividersbn} and Lemma \ref{toriofclassgrps} we obtain that
every element of order $r$ is contained in a maximal torus $T$ of
order $\frac{1}{(4,q^n-1)}(q^{n-1}+(-1)^k)(q+(-1)^k)$ and $k=n-1$ if
$k$ is odd. Then 2 does not divide $\vert T\vert$ if and only if
condition 2 or condition 3 of the proposition hold. If $\eta(k)=n$,
then $n$ is odd and, by \eqref{dividersbn} and Lemma
\ref{toriofclassgrps}, every element of order $r$ is contained in a
maximal torus $T$ of order $\frac{1}{(4,q^n-1)}(q^n-1)$. Clearly, $2$
does not divide $\vert T\vert$ if and only if $(4,q^n-1)=(q^n-1)_2$
and the proposition follows.
\end{proof}

\begin{prop}\label{adj2except}
Let $G$ be a finite simple exceptional group of Lie type over a
field of odd characteristic $p$. Let $r$ be an odd prime
divisor of $\vert G\vert$, $r\not=p$, and $k=e(r,q)$. Then $2$ and $r$ are
non-adjacent if and only if one of the following holds{\em:}
\begin{itemize}
\item[{\em 1.}] $G=G_2(q)$, $k\in\{3,6\}$.
\item[{\em 2.}] $G=F_4(q)$, $k=12$.
\item[{\em 3.}] $G=E_6(q)$, $k\in\{9,12\}$.
\item[{\em 4.}] $G={}^2E_6(q)$, $k\in\{12,18\}$.
\item[{\em 5.}] $G=E_7(q)$, and either $k\in\{7,9\}$ and $e(2,q)=2$, or
  $k\in\{14,18\}$ and~$e(2,q)=1$.
\item[{\em 6.}] $G=E_8(q)$, $k\in\{15,20,24,30\}$.
\item[{\em 7.}] $G={^3D_4(q)}$, $k=12$.
\item[{\em 8.}] $G={^2G_2(3^{2n+1})}$, and $r$
  divides $m_3(G,n)$ or~$m_4(G,n)$, where $m_i(G,n)$ are defined in
  Lemma~{\em\ref{SuzReeDivisors}}.
\end{itemize}
\end{prop}

\begin{proof}
The proposition is immediate from Lemmas~\ref{toriofexcptgrps}
and~\ref{SuzReeDivisors}.
\end{proof}

\section{On connection between the properties of the prime graph $GK(G)$ and the
structure of $G$}\label{influence}

Let $GK(G)$ be a prime graph of a finite group $G$. Obviously, the
spectrum $\omega(G)$ uniquely determines the structure of $GK(G)$.
Denote by $s(G)$ the number of connected components of $GK(G)$ and
by $\pi_i(G)$, $i=1,\dots, s(G)$, the $i$th connected component of
$GK(G)$. If $G$ has  even order then put $2\in\pi_1(G)$. Denote by
$\omega_i(G)$ the set of numbers $n\in\omega(G)$ such that each
prime divisor of $n$ belongs to $\pi_i(G)$.

Gruenberg and Kegel obtained the following description of finite groups with disconnected
prime graph.

\begin{ttt}\label{GKT}
{\em (Gruenberg~---~Kegel Theorem) (see \cite{Wil})}  If $G$ is a finite group with
disconnected graph $GK(G)$, then one of the following statements holds{\em:}

{\rm (a)}\ $s(G)=2$ and $G$ is a Frobenius group{\em;}

{\rm (b)}\ $s(G)=2$ and $G$ is a $2$-Frobenius group, i.e., $G=ABC$, where $A$, $AB$ are
normal subgroups of~$G$; $AB$, $BC$ are Frobenius groups with cores~$A$, $B$ and
complements~$B$, $C$ respectively{\em;}

{\rm (c)}\ there exists a nonabelian simple group~$S$ such that
$S\leq \overline{G}=G/K\leq \operatorname{Aut}(S)$, where $K$ is a
maximal normal soluble subgroup of $G$. Furthermore, $K$ and
$\overline{G}/S$ are $\pi_1(G)$-subgroups, the graph~$GK(S)$ is
disconnected, $s(S)\geqslant s(G)$, and for every~$i$, $2\leqslant
i\leqslant s(G)$, there is $j$, $2\leqslant j\leqslant s(S)$, such
that $\omega_i(G)=\omega_j(S)$.
\end{ttt}

Together with the classification of finite simple groups with
disconnected prime graph obtained by Williams and Kondrat'ev (see
\cite{Wil} and \cite{Kon}) the Gruenberg~---~Kegel theorem implied a
series of important corollaries (see, for example, \cite[Theorems
3--6]{Wil} and \cite[Theorems 2--3]{Kon}). In last years this
theorem is used for the proving of recognizability of finite groups by
spectrum (for details see \cite{Maz} and \cite{Vas}).

The proof of the Gruenberg~---~Kegel Theorem substantially used the
fact that in a group $G$ (if its order is even) there exists an
element of odd prime order disconnected in $GK(G)$ with a prime $2$.
It turns out that  disconnectedness of $GK(G)$ could be successfully
replaced in most cases by a ~weaker condition for the prime~$2$ to
be nonadjacent to at least one odd prime.

Denote by $t(G)$ the maximal number of prime dividers of the order of
$G$ pairwise non-adjacent in
$GK(G)$. In other words, $t(G)$ is a maximal number of vertices in independent sets
of $GK(G)$ (a set of vertices is called {\em independent} if  its elements
are pairwise non-adjoint). In graph theory this number is usually called an
{\em independence number} of
the graph. By analogy we denote by $t(r,G)$ the maximal number of vertices in
independent sets of $GK(G)$ containing the prime $r$. We call this number an
{\em $r$-independence number}.

\begin{ttt}\label{motiv}
{\em (see \cite{Vas})}  Let $G$ be a finite group satisfying two
conditions{\em:}

{\rm (a)} there exist three primes in $\pi(G)$ whose are pairwise
non-adjacent in $GK(G)$, i.~e.,~${t(G)\ge3}${\em;}

{\rm (b)} there exists an odd prime in $\pi(G)$ which is
non-adjacent to prime $2$ in $GK(G)$, i.~e., ${t(2,G)\ge2}$.

Then there exists a finite nonabelian simple group $S$ such that
$S\leq \overline{G}=G/K\leq \operatorname{Aut}(S)$ for maximal
normal soluble subgroup~$K$ of $G$. Furthermore, $t(S)\geqslant
t(G)-1$ and one of the following statements holds:

{\rm (1)} $S\simeq Alt_7$ or $A_1(q)$ for some odd $q$ and~$t(S)=t(2,S)=3$.

{\rm (2)} For every prime $r\in\pi(G)$ non-adjacent in $GK(G)$ to $2$ the Sylow
$r$-subgroup of $G$ is isomorphic to the Sylow $r$-subgroup of $S$. In particular,
$t(2,S)\ge t(2,G)$.
\end{ttt}

Note that the condition (a) of Theorem \ref{motiv} easily implies the
insolubility of
group $G$ and so, by the Feit~---~Thompson Odd Theorem, implies that $G$ is of even
order. Moreover, the condition (a) can be replaced by the weaker condition of the
insolubility of $G$ (see \cite[Propositions 2--3]{Vas}). Together with \cite[Theorem
2]{LM} Theorem \ref{motiv} implies the following statement.

\begin{ttt}\label{twoindepdescrib}
{\em (see \cite[Proposition 4]{Vas})} Let $L$ be a finite simple group with
$t(2,L)\ge 2$ non-isomorphic to $A_2(3)$, ${}^2A_3(3)$, $C_2(3)$ and $Alt_{10}$.
Let $G$ be a finite group with $\omega(G)=\omega(L)$. Then for a group $G$ the
conclusion of Theorem {\em \ref{motiv}} holds true. In particular, $G$ has a
unique
nonabelian composition factor.
\end{ttt}

As a matter of fact, these recent results (Theorem \ref{motiv} and
Theorem \ref{twoindepdescrib}) gave the motivation for the present
work. In Section \ref{sporadicalt} using results of Sections
\ref{prlmn}--\ref{two} we calculate the values of the independence
numbers and the 2-independence numbers for all finite nonabelian
simple groups. Furthermore, for every finite nonabelian simple group
of Lie type over the field of characteristic $p$ we also determine
$p$-independence number $t(p,G)$.

\section{Independence numbers}\label{sporadicalt}

For a finite group $G$ denote by $\rho(G)$ (by $\rho(r,G)$) some independence set in
$GK(G)$ (containing $r$) with maximal number of vertices. Thus, $|\rho(G)|=t(G)$ and
$|\rho(r,G)|=t(r,G)$. For a given finite nonabelian simple group  we point out some
$\rho(G)$ and $\rho(2,G)$ (obviously, they may be not uniquely defined) and so
determine $t(G)$ and~$t(2,G)$.

\begin{prop}\label{indsporadic}
Let $G$ be a sporadic group. Then $\rho(G)$, $t(G)$, $\rho(2,G)$ and
$t(2,G)$ are listed in Table~{\em\ref{sporadic}}. Furthermore,
$\rho(2,G)$ is uniquely determined.
\end{prop}

\begin{proof}
The results are easy to obtain using \cite{Atlas} or \cite{GAP}.
\end{proof}

Now we deal with simple alternating groups.

\begin{prop}\label{indalt}
Let $G=Alt_n$ be an alternating group of degree $n\geqslant 5$. Let
$s'_n$ be the largest prime which does not exceed $n/2$ and $s''_n$
be the smallest prime greater than $n/2$. Denote by $\tau(n)$ the
set $\{s~|~s \mbox{ is a prime, }n/2 < s\leqslant n\}$ and by
$\tau(2,n)$ the set $\{s~|~s \mbox{ is a prime, }n-3\leqslant
s\leqslant n\}$. Then $\rho(G)$, $t(G)$, $\rho(2,G)$ and $t(2,G)$
are listed in Table~{\em\ref{twoindalt}}. Furthermore, $\rho(2,G)$
is uniquely determined.
\end{prop}

\begin{proof}
The result is immediate corollary of Proposition \ref{adjalt}.
\end{proof}

Now we consider groups of Lie type. Since the problem here is more
complicated we divide the process of its solution into several natural steps.
The
main tools will be a Zsigmondy Theorem, Lemma~\ref{SuzReeDivisors}
and our results from Sections~\ref{odds}-\ref{two}. Recall that for
a given natural number $q$ we denote by $r_n$ some primitive prime
divisor of $q^n-1$ (if it exists). Note that a primitive prime
divisor $r_{2m}$ of $q^{2m}-1$ divides $q^m+1$ and does not divide
$q^k+1$ for every natural $k<m$. Since arguments in even characteristic of the
 field of definition are suitable for every characteristic $p$ we start with
determination of~${t(p,G)}$.

\begin{prop}\label{tpclassic}
Let $G$ be a finite simple group of Lie type over a field of
characteristic $p$ and order $q$. Then $\rho(p,G)$ and $t(p,G)$ are
listed in Table~{\em\ref{pindependenceclass}} for classical groups
and Table~{\em\ref{pindependenceex}} for exceptional groups.
\end{prop}

\begin{proof} (1) $G=A^{\varepsilon}_{n-1}(q)$.

Let $G=A_1(q)$. Since $A_1(2)$ and $A_1(3)$ are not simple we may
assume that $q>3$.  By Zsigmondy Theorem, if $q\not=2^t-1$, then there
exist primitive prime divisors $r_1$ and $r_2$ of $q-1$ and
$q+1$. They are non-adjacent in $GK(G)$. By Proposition
\ref{adjcharclass}, they are non-adjacent to $p$. Since $q>3$, for any
$q=2^t-1$ there exists odd primitive prime divisor $r_1$ of $q-1$
which is not adjacent to $2$ in $GK(G)$. Since in this case $e(2,q)=2$
and $q+1=2^t$, the prime $2$ is a unique primitive prime divisor of
$q^2-1$. Hence $\rho(p,G)=\{p,r_1,r_2\}$ in this case too.

Suppose that $n>2$. In order to consider groups $A_{n-1}(q)$ and ${}^2A_{n-1}(q)$
together we define new function:
$$\nu_{\varepsilon}(m)=\left\{
\begin{array}{rl}
m &\text{, if either }\varepsilon=+,\text{ or }\varepsilon=-\text{ and }m\equiv 0(\mod 4),\\
\frac{m}{2}& \text{, if }\varepsilon=-\text{ and }m\equiv 2(\mod 4),\\
2m&\text{, if }\varepsilon=-\text{ and }m\equiv1(\mod 2).\\
\end{array}\right.$$
Obviously, that $\nu_{\varepsilon}(m)$ is the identity function if
$\varepsilon=+$, and $\nu_{\varepsilon}(m)=\nu(m)$ if
$\varepsilon=-$. It is easy to check that $\nu_{\varepsilon}$ is a
bijection on $\mathbb{N}$, thus $\nu^{-1}_{\varepsilon}$ is well
defined.

Let $n=3$ and $G=A^{\varepsilon}_2(q)$. Then, by Proposition
\ref{adjcharclass}, it follows that $r_{\nu_\varepsilon^{-1}(2)}$ and
$r_{\nu_\varepsilon^{-1}(3)}$ are non-adjacent to $p$. Moreover, if
$(q-\varepsilon1)_3=3$, then, by Proposition \ref{adjcharclass}, we
have that 3 is non-adjacent to $p$ and, by Propositions
\ref{bigprimesan} and \ref{bigprimesantw}, in this case 3 is
non-adjacent to $r_{\nu_\varepsilon^{-1}(2)}$ and
$r_{\nu_\varepsilon^{-1}(3)}$. By Zsigmondy Theorem
$r_{\nu_\varepsilon^{-1}(3)}$ exists for all $q$, except
$\varepsilon=-$ and $q=2$. But ${}^2A_2(2)$ is not simple and we do
not consider this group. Therefore $r_{\nu_\varepsilon^{-1}(3)}$
exists for all simple groups in this case. Again by Zsigmondy Theorem,
$r_{\nu_{\varepsilon^{-1}(2)}}$ exists if and only if $q+\varepsilon1$
is not a power of $2$, i.e., it is not a Mersenne prime in case
$\varepsilon=+$, and it is not a Fermat prime or $9$ in case
$\varepsilon=-$. Thus, we have
$$\rho(p,G)=\left\{\begin{array}{rl}
\{p,3,r_{\nu_{\varepsilon}^{-1}(2)},r_{\nu_{\varepsilon}^{-1}(3)}\},
\text{ if }(q-\varepsilon1)_3=3\text{ and }q+\varepsilon\neq 2^k;\\
\{p,r_{\nu_{\varepsilon}^{-1}(2)},r_{\nu_{\varepsilon}^{-1}(3)}\}, \text{
if }(q-\varepsilon1)_3\neq 3\text{ and }q+\varepsilon\neq 2^k;\\
\{p,3,r_{\nu_{\varepsilon}^{-1}(3)}\}, \text{ if }(q-\varepsilon1)_3=3\text{
and }q+\varepsilon=2^k;\\
\{p,r_{\nu_{\varepsilon}^{-1}(3)}\}, \text{ if }(q-\varepsilon1)_3\neq 3\text{ and
}q+\varepsilon=2^k.
\end{array}\right.$$
Note that we avoid to use the functions $e(r,q)$ and
$\nu_{\varepsilon}$ in the tables in Section~\ref{tables}.

Let $G=A_4(2), A_5(2)\text{ or }{}^2A_3(2)$. Since there are no primitive prime
divisors of ${2^6-1}$, by Proposition \ref{adjcharclass} we have that
$\rho(2,A_5(2))=\{2,31\}$, $\rho(2,A_6(2))=\{2,127\}$, and
$\rho(2,{}^2A_3(2))=\{2,5\}.$

In all other cases by Zsigmondy Theorem there exist primitive prime
divisors $r_{\nu_{\varepsilon}^{-1}(n-1)}$
and~$r_{\nu_{\varepsilon}^{-1}(n)}$. By Lemma \ref{toriofclassgrps}
we have $r_{\nu_{\varepsilon}^{-1}(n-1)},
r_{\nu_{\varepsilon}^{-1}(n)}\in\pi(G)$. By Propositions
\ref{adjanspl} and \ref{adjantwst} they are non-adjacent in $GK(G)$.
Now Proposition \ref{adjcharclass} yields that for
$G=A^{\varepsilon}_{n-1}(q)$ we have
$\rho(p,G)=\{p,r_{\nu_{\varepsilon}^{-1}(n-1)},
r_{\nu_{\varepsilon}^{-1}(n)}\}$.

(2) $G=C_n(q)$ or $B_n(q)$.

In view of Proposition \ref{adjbn}, Proposition \ref{adjcharclass}, and
Proposition \ref{adj2bn}, we have that the prime graphs of $C_n(q)$ and $B_n(q)$
coincide. So we consider these groups together and, for brevity, use the symbol
$C_n(q)$ in both cases.

Let $G=C_3(2)$. Since there are no primes $r$ with $e(r,2)=6$, only
$7$ as the primitive prime divisor of $2^3-1$ is not adjacent to
$2$.

Let $G=C_n(q)$, $n\ge 2$ and $(n,q)\neq(3,2)$. If $n$ is even, by
Proposition \ref{adjcharclass} only primitive prime divisors of
$q^n+1$ are non-adjacent to $p$. Thus, in this case
$\rho(G)=\{p,r_{2n}\}$. If $n$ is odd, the Propositions
\ref{adjbn} and \ref{adjcharclass} yield
$\rho(G)=\{p,r_n,r_{2n}\}$.

(3) $G=D^{\varepsilon}_n(q)$.

With respect to well-known isomorphisms of groups of small Lie rank,
we may suppose that $n\ge 4$. Let $n=4$ and $q=2$. Since there are
no primitive prime divisors of $2^6-1$ we have that only $7$ as a
primitive prime divisor of $2^3-1$ is non-adjacent to $2$ in case
$G=D_4(2)$ and only $7$ and $17$ as primitive divisors of $2^3-1$
and $2^8-1$ are non-adjacent to $2$ in case $G={}^2D_4(2)$. All
others possibilities could be easily described in the direct
accordance to Proposition~\ref{adjdn} and
Proposition~\ref{adjcharclass}. The results of such consideration one can
see in Table~\ref{pindependenceclass}.

(4) The result for exceptional groups distinct from Suzuki and Ree
groups can be obtained by using Proposition~\ref{adjexcept},
Proposition~\ref{adjcharexcept} and the Zsigmondy Theorem.

(5) $G$ is a finite simple Suzuki or Ree group.

Let $G={}^2B_2(2^{2m+1})$, and $s_i$ be a prime divisor of
$m_i(B,n)$, where $i=1,2,3$ (see Lemma~\ref{SuzReeDivisors}). By
Proposition~\ref{adjsuzree} primes $s_i$ and $s_j$ are adjacent if
and only if $i=j$. On the other hand, every $s_i$, where $i=1,2,3$,
is non-adjacent to~$p=2$. Thus, $\rho(p,G)=\{p,s_1,s_2,s_3\}$ and
$t(p,G)=4$ in this case. The same arguments one can apply to Ree
groups and prime divisors of $m_i(G,n)$ and $m_i(F,n)$.
\end{proof}

In general, a primitive prime divisor $r_m$ of $q^m-1$ could be chosen
by several ways. Thus, the set $\rho(p,G)$ could be not uniquely
determined for a finite simple group $G$ of Lie type.  However, it
turns out that the values of $e(r,q)$ for all primitive prime divisors
$r$ from $\rho(p,G)$ are invariants for a given group $G$ of Lie type.

\begin{prop}\label{uniquepind}
Let $G$ be a finite simple group of Lie type over a field of
characteristic $p$ and order $q$. Assume that $G$ is not isomorphic
to ${}^2B_2(q)$, ${}^2G_2(q)$, ${}^2F_4(2)'$, and ${}^2F_4(q)$. Let
$\rho(p,G)=\{p,s_1,s_2,\ldots,s_m\}$ be an independence set in
$GK(G)$ containing $p$ with maximal number of vertices
and~$k_i=e(s_i,q)$. Then the set $\{k_1,\ldots,k_m\}$ is uniquely
determined.
\end{prop}

\begin{proof}
It follows from the results of Sections \ref{odds}-\ref{two}.
\end{proof}

If we change primitive prime divisors on divisors of numbers defined
in Lemma~\ref{SuzReeDivisors}, we obtain a similar statement for
Suzuki and Ree groups.

\begin{prop}\label{uniquepindsuzree}
Let $G$ be a finite simple Suzuki or Ree group over a field of
characteristic $p$. Let~$\rho(p,G)=\{p,s_1,\ldots,s_k\}$.
\begin{itemize}
\item[{\em 1.}] If $G={^2B_2(2^{2n+1})}$, then, up to reodering, $s_i$
  divides $m_i(B,n)$. In particular,~$k=3$.
\item[{\em 2.}] If  $G={^2G_2(3^{2n+1})}$, then, all of $s_i$ are odd
  and, up to reodering, $s_i$ divides $m_i(G,n)$. In particular,~$k=4$.
\item[{\em 3.}] If $G={^2F_4(2^{2n+1})}$, then, up to reodering, $s_i$
  divides $m_{i+2}(F,n)$ and $s_1\not=3$. In particular,~$k=3$.
\end{itemize}
Numbers $m_i(B,n)$, $m_i(G,n)$, and $m_i(F,n)$ are defined in
Lemma~{\em\ref{SuzReeDivisors}}.
\end{prop}

Now we determine $t(2,G)$. Obviously, for a group of Lie type over a
field of even characteristic $p$-independence number and
$2$-independence number coincide. Thus, we may assume that $G$ is
defined over the field of odd characteristic.

\begin{prop}\label{t2classic}
 Let $G$ be a finite simple group of Lie type
over the field of odd characteristic $p$ and order $q$. Then
$\rho(2,G)$ and $t(2,G)$ are listed in Table~{\em
\ref{2independenceclass}} for classical groups and in Table~{\em
\ref{2independenceex}} for exceptional groups.
\end{prop}

\begin{proof} (1) $G=A^{\varepsilon}_{n-1}(q)$.

Let $G=A_1(q)$ and $q>3$. Since $q$ is odd, every prime divisor
$r\neq p$ of $|G|$ divides $(q-1)/2$ or $(q+1)/2$. If $e(2,q)=1$,
then $2$ is non-adjacent to some prime divisor $r_2$ of $(q+1)/2$
and $\tau(G)=\{2,p,r_2\}$. If $e(2,q)=2$, then $2$ is non-adjacent
to some divisor $r_1$ of $(q-1)/2$ and $\tau(G)=\{2,p,r_1\}$.

Let $G=A^{\varepsilon}_{n-1}(q)$ and $n\ge 3$. If $(q-\varepsilon1)_2<n_2$, then by
Propositions \ref{bigprimesan} and \ref{bigprimesantw} only primitive prime divisor
$r_{\nu_{\varepsilon}^{-1}(n-1)}$ is non-adjacent to $2$. Note that in this case
$n_2\ge4$, since for $n_2\le2$ the inequality $(q-\varepsilon 1)_2<n_2$ is
impossible. By Zsigmondy Theorem primitive prime divisor
$r_{\nu_{\varepsilon}^{-1}(n-1)}$ always exists.
Thus,~$\rho(2,G)=\{2,r_{\nu_{\varepsilon}^{-1}(n-1)}\}$.

If $(q-\varepsilon1)_2>n_2$ or $(q-\varepsilon1)_2=n_2=2$, then by Propositions
\ref{bigprimesan} and \ref{bigprimesantw}  every primitive prime divisor
$r_{\nu_{\varepsilon}^{-1}(n)}$ is non-adjacent to $2$. Therefore,
$\rho(2,G)=\{2,r_{\nu_{\varepsilon}^{-1}(n)}\}$.

At last, let $(q-\varepsilon1)_2=n_2>2$. By Propositions \ref{bigprimesan} and
\ref{bigprimesantw} only primitive prime divisors $r_{\nu_{\varepsilon}^{-1}(n-1)}$
and $r_{\nu_{\varepsilon}^{-1}(n)}$ are non-adjacent to $2$. On the other hand, by
Propositions \ref{adjanspl} and \ref{adjantwst}, primes
$r_{\nu_{\varepsilon}^{-1}(n-1)}$ and $r_{\nu_{\varepsilon}^{-1}(n)}$ are
non-adjacent. Thus,
$\rho(2,G)=\{2,r_{\nu_{\varepsilon}^{-1}(n-1)},r_{\nu_{\varepsilon}^{-1}(n)}\}$.

(2) $G=C_n(q)$ or $B_n(q)$.

The results from the table are directly following from Proposition
\ref{adj2bn}.

(3) $G=D^{\varepsilon}_n(q)$.

The results again are the direct corollary of Proposition \ref{adj2dn}. Note that
the equality $t(2,G)=2$ is true for most groups of type $D_n$ over fields of odd
order. Exceptions are following: $n$ is odd and $q\equiv 5(\mod 8)$ for $G=D_n(q)$,
and $q\equiv 3(\mod 8)$ for $G={}^2D_n(q)$.

(4) Now to complete the proof of the proposition one can use
Proposition \ref{adj2except} and Lemma~\ref{SuzReeDivisors} for
Suzuki and Ree groups, and the Zsigmondy Theorem for other
exceptional groups.
\end{proof}

Now we consider the uniqueness of $\rho(2,G)$. The situation here is
very similar to the situation with $\rho(p,G)$.

\begin{prop}\label{uniquetwoind}
Let $G$ be a finite simple group of Lie type over a field of odd
characteristic $p$ and order $q$. Assume that $G$ is not isomorphic
to $A_1(q)$, and ${}^2G_2(q)$. Let
$\rho(2,G)=\{2,s_1,s_2,\ldots,s_m\}$ be an independence set in
$GK(G)$ containing $2$ with maximal number of vertices
and~$k_i=e(s_i,q)$. Then the set $\{k_1,\ldots,k_m\}$ is uniquely
determined.
\end{prop}

\begin{proof}
It follows from the results of Sections \ref{odds}--\ref{two}.
\end{proof}

\begin{prop}\label{uniquetwoindsuzree}
Let $G$ be either $A_1(q)$ with $q$ odd, or
  ${}^2G_2(3^{2n+1})$ and let $\rho(2,G)=\{2,s_1,s_2\}$.
\begin{itemize}
\item[{\em 1.}] If $G=A_1(q)$, $q$ odd, then, up to renumbering,
  $s_1=p$ and~$e(s_2,q)=3-e(2,q)$.
\item[{\em 2.}] If $G={^2G_2(3^{2n+1})}$, then, up to
  renumbering,~$s_i=m_{i+2}(G,n)$.
\end{itemize}
Numbers $m_i(G,n)$ are defined in Lemma~{\em\ref{SuzReeDivisors}}.
\end{prop}

\begin{proof}
It follows from the results of Sections \ref{odds}--\ref{two}.
\end{proof}

Our last task is to determine for every finite simple group $G$ of
Lie type some independent set $\rho(G)$ in $GK(G)$ with maximal
number of vertices.

\begin{prop}\label{tan}
Let $G$ be a finite simple group of Lie type over the field of
characteristic $p$ and order $q$. Then $\rho(G)$ and $t(G)$ are
listed in Table~{\em\ref{independence}} for classical groups and in
Table~{\em\ref{indepexcept}} for exceptional groups.
\end{prop}

\begin{proof} (1) $G=A^{\varepsilon}_{n-1}(q)$.

If $G=A_1(q)$ or $G=A^{\varepsilon}_2(q)$, then arguing as in proof of
Proposition \ref{tpclassic} we obtain $\rho(G)=\rho(p,G)$.

Suppose $n=4$. In view of Proposition \ref{tpclassic} and Table
\ref{pindependenceclass} we have that $t(p,G)=3$ in all cases,
except ${}^2A_3(2)$. By Proposition \ref{t2classic} and Table
\ref{2independenceclass} we obtain that $t(2,G)\le3$. By
Propositions \ref{bigprimesan} and \ref{bigprimesantw} it follows
that $t(r_{\nu_\varepsilon^{-1}(1)},G)\le 3$. Furthermore, there exist
at most three other primitive divisors
$r_{\nu_\varepsilon^{-1}(2)}$, $r_{\nu_\varepsilon^{-1}(3)}$, and
$r_{\nu_\varepsilon^{-1}(4)}=r_4$. So $t(G)= t(p,G)=3$ except the
case: $t({}^2A_3(2))=t(2,{}^2A_3(2))=2$.

Let $G=A^{\varepsilon}_{n-1}(q)$, $n\ge 5$ and $q\neq 2$. By
Zsigmondy Theorem there exist primitive prime divisors
$r_{\nu^{-1}_{\varepsilon}(k)}$ for every $k>2$. Denote by $m$ the
number $[n/2]$, i.e., the integral part of $n/2$. By Propositions
\ref{adjanspl} and \ref{adjantwst} the set
$$\rho=\{r_{\nu^{-1}_{\varepsilon}(m+1)},r_{\nu^{-1}_{\varepsilon}(m+2)},
\ldots,r_{\nu^{-1}_{\varepsilon}(n)}\}$$ is an independent set in
$GK(G)$. On the other hand, by Propositions \ref{adjanspl},
\ref{adjantwst}, \ref{adjcharclass}, \ref{bigprimesan} and
\ref{bigprimesantw}, every two prime divisors from
$\pi\left(q\prod_{i=1}^m(q^i-(\varepsilon1)^i)\right)$ are adjacent
in $GK(G)$. Furthermore, each of them is adjacent to at least one
number from $\rho$. Since every prime $s\in \pi(G)\setminus
\left(\rho\cup\pi\left(q\prod_{i=1}
^m\left(q^i-(\varepsilon1)^i\right)\right)\right)$ is of
the form $r_i$ for some $i>m$, it follows that every independent set
in $GK(G)$ with prime divisor from
$\pi\left(q\prod_{i=1}^m(q^i-1)\right)$ contains at most
$\vert\rho\vert$ vertices. Thus, $\rho(G)=\rho$
and~$t(G)=\vert\rho\vert=[(n+1)/2]$.

Let $q=2$. Since $\nu^{-1}_{+}(6)=\nu^{-1}_{-}(3)=6$, the results of
previous paragraph hold true for $A_{n-1}(2)$ with $n\ge 12$ and
${}^2A_{n-1}(2)$ with $n\ge 6$. If $n=5$ and $\varepsilon=-$, then
for  ${}^2A_4(2)$ one can check that $\rho(G)=\rho(2,G)=\{2,5,11\}$
and $t(G)=3$. So we can suppose that $G=A_{n-1}(2)$. If $n=5,6$ then
$\rho(G)=\{r_3,r_4,r_5\}=\{5,7,31\}$ and $t(G)=3$. If $7\leqslant n
\leqslant 11$ then we need to eliminate divisors of $2^6-1$ from
$\rho(G)$ and at rest argue like in previous paragraph. In this case
$\rho(G)=\{r_i~|~i\neq 6, [n/2]<i\leqslant n\}$ and
$t(G)=[(n-1)/2]$.

(2) $G=C_n(q)$ or $B_n(q)$. Recall that $GK(C_n(q))=GK(B_n(q))$ and we
    consider these groups together.

$G=C_2(q)$, $q>2$. Since every two prime divisors of $q(q^2-1)$ are
adjacent, we obtain $\rho(G)=\rho(p,G)=\{p,r_4\}$.

Let $n\geqslant 3$ if $q>2$, and $n\ge 7$ if $q=2$. Define the set
$\rho$ as follows:
$$\rho=\{r_{2i}~|~[n+1/2]< i\le n\}\cup\{r_{i}~|~[n/2]<i\le
n,i\equiv{1}(\mod 2)\}.$$ Using results of Sections
\ref{prlmn}--\ref{two} and arguments as in proof for groups
$A^{\varepsilon}_{n-1}(q)$ we obtain that $\rho(G)=\rho$ and as easy
to verify $t(G)=[(3n+5)/4]$.

If $q=2$ and $3\le n\le 5$, all arguments stand the same but we have
to eliminate the divisors of type $r_6$ from $\rho$. Thus, in this
case $t(G)=[(3n+1)/4]$. At last, if $(n,q)=(6,2)$ we have to
eliminate a divisor of type $r_6$ from $\rho$, but instead it we can
add a primitive prime divisor $7$ of $2^3-1$. Thus,
$t(G)=[(3n+5)/4]$ as in the common situation.

(3) $G=D^{\varepsilon}_n(q)$.

We argue like in two previous parts of proof. Using Proposition
\ref{adjdn} we obtain a set $\rho(G)$ in common situation and then
consider some exceptions arising with respect to exceptions in the
Zsigmondy Theorem.

(4) We obtain results for exceptional groups by using
Propositions~\ref{adjexcept} and~\ref{adjsuzree}, Tables
\ref{pindependenceex} and \ref{2independenceclass}, and the
Zsigmondy Theorem and Lemma \ref{SuzReeDivisors}.
\end{proof}

\section{Applications}\label{applications}

Here we apply our results in a spirit of Section \ref{influence}.
First of all our investigations show that the condition $t(2,G)>1$
is realized for an extremely wide class of finite simple groups.

\begin{ttt}\label{theorem}
Let $G$ be a finite nonabelian simple group with $t(2,G)=1$, then
$G$ is an alternating group $Alt_n$ with $\tau(2,n)=\{s~|~s \mbox{
is a prime, }n-3\leqslant s\leqslant n\}=\varnothing$.
\end{ttt}

\begin{proof} See Propositions \ref{indsporadic}--\ref{tan} and corresponding tables in Section~\ref{tables}.
\end{proof}

Thus, we may apply results of Theorems \ref{motiv} and
\ref{twoindepdescrib} as follows.

\begin{cor}\label{describ}
Let $L$ be a finite nonabelian simple group distinct from $A_2(3)$,
${}^2A_3(3)$, $C_2(3)$, $Alt_{10}$ and $Alt_n$ with
$\tau(2,n)=\varnothing$. Let $G$ be a finite group with
$\omega(G)=\omega(L)$. Then there exists a finite nonabelian simple
group $S$ such that $S\leq \overline{G}=G/K\leq
\operatorname{Aut}(S)$ for maximal normal soluble subgroup~$K$ of
$G$. Furthermore, $t(S)\geqslant t(G)-1$ and one of the following
statements holds{\em:}

{\rm (1)} $S\simeq Alt_7$ or $A_1(q)$ for some odd $q$ and
$t(S)=t(2,S)=3$.

{\rm (2)} For every prime $r\in\pi(G)$ non-adjacent in $GK(G)$ to
$2$ the Sylow $r$-subgroup of $G$ is isomorphic to the Sylow
$r$-subgroup of $S$. In particular, $t(2,S)\geqslant t(2,G)$.
\end{cor}

\begin{proof}
It is direct corollary of Theorem \ref{twoindepdescrib} and Theorem
\ref{theorem}.
\end{proof}

If $G$ is one of the groups $A_2(3)$, ${}^2A_3(3)$, $C_2(3)$,
$Alt_{10}$, then $t(2,G)\ge2$ and we also have

\begin{cor}\label{mazquest}
Let $L$ be a finite simple group distinct from $Alt_n$ with
$\tau(2,n)=\varnothing$. Let $G$ be a finite group with
$\omega(G)=\omega(L)$. Then $G$ has at most one nonabelian
composition factor.
\end{cor}

For an arbitrary subset $\omega$ of the set $\mathbb{N}$ of natural
numbers denote by $h(\omega)$ the number of pairwise nonisomorphic
finite groups $G$ such that $\omega(G)=\omega$. We say that for a
finite group $G$ {\em the recognition problem is solved} if we know
the value of $h(\omega(G))$ (for brevity, $h(G)$). In particular,
the group $G$ is said to be {\it recognizable by spectrum} (briefly,
recognizable) if $h(G)=1$. In last twenty years the recognition
problem is solved for many finite nonabelian simple and almost
simple groups (for the review of results in this field see
\cite{Maz}). But most of those groups have the disconnected prime
graph, since the Gruenberg~---~Kegel Theorem and classification of
finite simple groups with disconnected prime graph obtained by
Williams and Kondrat'ev were an important part of proof. The main
theorem of \cite{Vas} (Theorem \ref{motiv} here) and the results of
the present work give a possibility to deal with a groups whose
prime graph is connected. A finite nonabelian simple group $L$ is
said to be {\it quasirecognizable by spectrum} if every finite group
with the same spectrum has exactly one nonabelian composition factor
$S$ and $S$ is isomorphic to $L$. So the investigations of
quasirecognizability is an important step in the determination,
whether the given group is recognizable by spectrum. In \cite{Vas}
the first author gave a sketch of a proof for the following
statement.

\begin{ttt}\label{quasi} {\rm(\cite[Proposition 5]{Vas})}
Let $L={}^2D_n(q)$, $q=2^k$, $k,n$ are natural numbers, $n$ is
even and $n\geqslant 16$. Then $L$ is quasirecognizable by spectrum.
\end{ttt}

Actually, this result was proven modulo the results of Section
\ref{sporadicalt} of the present article. Thus, now it is completely
proven.

Now we would like to emphasize one statement that we have mentioned
before.

\begin{prop}\label{coincide}
Let $G=B_n(q)$ and $H=C_n(q)$. Then the prime graphs $GK(G)$ and
$GK(H)$ coincide.
\end{prop}

\begin{proof}
It follows from the results of Sections~\ref{odds}--\ref{two}.
\end{proof}

At last, we mention one recent result on prime graphs of finite
groups. Recall that a set of vertices of a graph is said to be a
clique if all vertices from this set are pairwise adjacent. In
\cite[Theorem 1]{LM} Lucido and Moghaddamfar describe all finite
nonabelian simple groups whose prime graph connected components are
cliques. We check their list of such groups using results of the
present article. Unfortunately, there are some mistakes in the list.
As a matter of fact, for groups $G=A^{\varepsilon}_2(q)$, where
$q=2^k-\varepsilon{1}$ and $(q-\varepsilon{1})_3=3$, we have
$3,p\in\pi_1(G)$ and $3p\not\in\omega(G)$. So the component
$\pi_1(G)$ is not a clique for those groups contrary to the
statement of Theorem~1 in~\cite{LM}. Below in Corollary
\ref{cliques} we give a revised list of such groups.

\begin{cor}\label{cliques}
Let $G$ be a finite nonabelian simple group and all connected
components of its prime graph $GK(G)$ are cliques. Then $G$ is one
of the following groups{\em:}
\begin{itemize}
\item[{\em 1.}] Sporadic groups $M_{11}$, $M_{22}$, $J_1$, $J_2$, $J_3$,
$HiS$.
\item[{\em 2.}] Alternating groups $Alt_n$, where $n=5, 6, 7, 9, 12, 13$.
\item[{\em 3.}] Groups of Lie type $A_1(q)$, where $q>3${\em;} $A_2(4)${\em;}
$A_2(q)$, where $(q-1)_3\neq 3$, $q+1=2^k${\em;}
${}^2A_3(3)${\em;} ${}^2A_5(2)$; ${}^2A_2(q)$, where $(q+1)_3\neq 3$,
$q-1=2^k${\em;} $C_3(2)$, $C_2(q)$, where $q>2${\em;} $D_4(2)${\em;}
${}^3D_4(2)${\em;}
${}^2B_2(q)$, where $q=2^{2k+1}${\em;} $G_2(q)$, where $q=3^k$.
\end{itemize}
\end{cor}

\begin{proof}
All connected component of prime graph $GK(G)$ of a finite group $G$
are cliques if and only if the number $s(G)$ of the components is
equal to the independence number $t(G)$ of $GK(G)$. For every finite
nonabelian simple group the values of $s(G)$ and $t(G)$ are now
known. Thus, using Tables 2a--2c in \cite{Maz} for the values of
$s(G)$ and Tables \ref{sporadic}--\ref{indepexcept} in Section
\ref{tables} of the present article for $t(G)$, we obtain the result
of the corollary.
\end{proof}
\newpage

\section{Resulting Tables}\label{tables}

In the tables below $n,k$ are assumed to be naturals. By
$[x]$ we denote the integral part of $x$. For a finite group $G$ we
denote by $\rho(G)$ (by $\rho(r,G)$) some independence set in
$GK(G)$ (containing $r$) with maximal number of vertices and put
$t(G)=|\rho(G)|$, $t(r,G)=|\rho(r,G)|$. In Table \ref{twoindalt} by
$\tau(n)$ we denote the set $\{s~|~s \mbox{ is a prime, }n/2 <
s\leqslant n\}$ and by $\tau(2,n)$ we denote the set $\{s~|~s \mbox{
is a prime, }n-3\leqslant s\leqslant n\}$. We denote $s'_n$ to be
the largest prime which does not exceed $n/2$ and $s''_n$ to be the
smallest prime greater than $n/2$. In Tables
\ref{pindependenceclass}--\ref{indepexcept} we assume $G$ to be a finite
nonabelian simple group of Lie type over a field of characteristic
$p$ and order $q$. By $r_m$ we define the primitive prime divisor of
$q^m-1$. If $p$ is odd then we say that $2$ is a primitive prime
divisor of $q-1$ if $q\equiv1(\mod4)$ and that $2$ is a primitive
prime divisor of $q^2-1$ if~$q\equiv -1(\mod4)$.

%\newpage

\begin{tab}\label{sporadic}{\bfseries Sporadic groups}\vspace{1\baselineskip}

{\small
\begin{tabular}{|r|c|l|c|l|}
 \hline $G$ & $t(G)$ & $\rho(G)$ & $t(2,G)$ & $\rho(2,G)$ \\
 \hline $M_{11}$ & $3$ & $\{3,5,11\}$ & $3$ & $\{2,5,11\}$ \\
 $M_{12}$ & $3$ & $\{3,5,11\}$ & $2$ & $\{2,11\}$ \\
 $M_{22}$ & $4$ & $\{3,5,7,11\}$ & $4$ & $\{2,5,7,11\}$ \\
 $M_{23}$ & $4$ & $\{3,7,11,23\}$ & $4$ & $\{2,5,11,23\}$ \\
 $M_{24}$ & $4$ & $\{5,7,11,23\}$ & $3$ & $\{2,11,23\}$ \\
 $J_{1}$ & $4$ & $\{5,7,11,19\}$ & $4$ & $\{2,7,11,19\}$ \\
 $J_{2}$ & $2$ & $\{5,7\}$ & $2$ & $\{2,7\}$ \\
 $J_{3}$ & $3$ & $\{5,17,19\}$ & $3$ & $\{2,17,19\}$ \\
 $J_{4}$ & $7$ & $\{7,11,23,29,31,37,43\}$ & $6$ & $\{2,23,29,31,37,43\}$ \\
 $\operatorname{Ru}$ & $4$ & $\{5,7,13,29\}$ & $2$ & $\{2,29\}$ \\
 $\operatorname{He}$ & $3$ & $\{5,7,17\}$ & $2$ & $\{2,17\}$ \\
 $\operatorname{McL}$ & $3$ & $\{5,7,11\}$ & $2$ & $\{2,11\}$ \\
 $\operatorname{HN}$ & $3$ & $\{7,11,19\}$ & $2$ & $\{2,19\}$ \\
 $\operatorname{HiS}$ & $3$ & $\{5,7,11\}$ & $3$ & $\{2,7,11\}$ \\
 $\operatorname{Suz}$ & $4$ & $\{5,7,11,13\}$ & $3$ & $\{2,11,13\}$ \\
 $\operatorname{Co}_{1}$ & $4$ & $\{7,11,13,23\}$ & $2$ & $\{2,23\}$ \\
 $\operatorname{Co}_{2}$ & $4$ & $\{5,7,11,23\}$ & $3$ & $\{2,11,23\}$ \\
 $\operatorname{Co}_{3}$ & $4$ & $\{5,7,11,23\}$ & $2$ & $\{2,23\}$ \\
 $\operatorname{Fi}_{22}$ & $4$ & $\{5,7,11,13\}$ & $2$ & $\{2,13\}$ \\
 $\operatorname{Fi}_{23}$ & $5$ & $\{7,11,13,17,23\}$ & $3$ & $\{2,17,23\}$ \\
 $\operatorname{Fi}'_{24}$ & $6$ & $\{7,11,13,17,23,29\}$ & $4$ & $\{2,17,23,29\}$ \\
 $\operatorname{O'N}$ & $5$ & $\{5,7,11,19,31\}$ & $4$ & $\{2,11,19,31\}$ \\
 $\operatorname{LyS}$ & $6$ & $\{5,7,11,31,37,67\}$ & $4$ & $\{2,31,37,67\}$ \\
 $F_{1}$ & $11$ & $\{11,13,17,19,23,29,31,41,47,59,71\}$ & $5$ & $\{2,29,41,59,71\}$ \\
 $F_{2}$ & $8$ & $\{7,11,13,17,19,23,31,47\}$ & $3$ & $\{2,31,47\}$ \\
 $F_{3}$ & $5$ & $\{5,7,13,19,31\}$ & $4$ & $\{2,13,19,31\}$ \\
   \hline
\end{tabular}}
\end{tab}

\newpage

\vfill
\begin{tab}\label{twoindalt}{\bfseries Simple alternating groups}\vspace{1\baselineskip}

{\small
\begin{tabular}{|c|l|c|c|c|c|}\hline
$G$&Conditions&$t(G)$&$\rho(G)$& $t(2,G)$&$\rho(2,G)$ \\ \hline
 $Alt_{n}$ & $n=5,6$ & $3$ & $\{2,3,5\}$ & $3$ & $\{2,3,5\}$\\
&$n=8$&$3$&$\{2,5,7\}$&$3$&$\{2,5,7\}$\\
& $n\ge 7$, $s'_n+s''_n > n$ & $|\tau(n)|+1$ & $\tau(n)\cup\{s'_n\}$ & $|\tau(2,n)|+1$ & $\tau(2,n)\cup\{2\}$\\
 & $n\ge 9$,  $s'_n+s''_n \le n$ & $|\tau(n)|$ & $\tau(n)$ & $|\tau(2,n)|+1$ & $\tau(2,n)\cup\{2\}$\\
 \hline
\end{tabular}}
\end{tab}
\vfill

\begin{tab}\label{pindependenceclass}{\bfseries $p$-independence
 numbers for finite simple classical groups}\vspace{1\baselineskip}

{\small
\begin{tabular}{|c|l|c|c|}\hline
$G$&Conditions&$t(p,G)$&$\rho(p,G)$\\ \hline
 $A_{n-1}(q)$ & $n=2$, $q>3$ & $3$ & $\{p,r_1,r_2\}$\\
 &$n=3$, $(q-1)_3=3$, and $q+1\not=2^k$&$4$&$\{p,3,r_2\not=2,r_3\}$\\
 &$n=3$, $(q-1)_3\not=3$, and $q+1\not=2^k$&$3$&$\{p,r_2\not=2,r_3\}$\\
 & $n=3$, $(q-1)_3=3\text{ and }q+1=2^k$ & $3$ & $\{p,3,r_3\}$\\
 & $n=3$, $(q-1)_3\neq 3\text{ and }q+1=2^k$ & $2$ & $\{p,r_3\}$\\
 & $n=6$, $q=2$& $2$ &$\{2,31\}$\\
 & $n=7$, $q=2$& $2$ &$\{2,127\}$\\
 & $n>3$ $(n,q)\neq(6,2),(7,2)$ & $3$ & $\{p,r_{n-1},r_n\}$ \\
 \hline
 ${}^2A_{n-1}(q)$ &$n=3$, $(q+1)_3=3$, and
$q-1\not=2^k$&$4$&$\{p,3,r_1\not=2,r_6\}$\\
&$n=3$, $(q+1)_3\not=3$, and $q-1\not=2^k$&$3$&$\{p,r_1\not=2,r_6\}$\\
 & $n=3$, $(q+1)_3=3\text{ and }q-1=2^k$ & $3$ & $\{p,3,r_6\}$\\
 & $n=3$, $(q+1)_3\neq 3\text{ and }q-1=2^k$ & $2$ & $\{p,r_6\}$\\
 & $n=4$, $q=2$ & $2$ & $\{2,5\}$ \\
 & $n\equiv 0(\mod{4})$, $(n,q)\neq(4,2)$  & $3$ & $\{p,r_{2n-2},r_n\}$ \\
 & $n\equiv 1(\mod{4})$  & $3$ &  $\{p,r_{n-1},r_{2n}\}$ \\
 & $n\equiv 2(\mod{4})$, $n\neq 2$ & $3$ & $\{p,r_{2n-2},r_{n/2}\}$ \\
 & $n\equiv 3(\mod{4})$, $n\neq 3$ & $3$ & $\{p,r_{(n-1)/2},r_{2n}\}$ \\
  \hline
 $B_n(q)$ or&$n=3$, $q=2$& $2$ & $\{2,7\}$\\
 $C_n(q)$& $n$ is even  & $2$ & $\{p,r_{2n}\}$\\
 & $n>1$ is odd, $(n,q)\neq(3,2)$ & $3$ & $\{p,r_n,r_{2n}\}$\\
 \hline
 $D_n(q)$ & $n=4$, $q=2$ & $2$ & $\{2,7\}$\\
 &$n\equiv0(\mod2)$, $n\ge 4$, $(n,q)\neq(4,2)$ & $3$ &
$\{p,r_{n-1},r_{2n-2}\}$\\
 &$n\equiv1(\mod2)$, $n> 4$ & $3$ &$\{p,r_{n},r_{2n-2}\}$\\
 \hline
 ${^2D_n(q)}$ & $n=4$, $q=2$ & $3$ & $\{2,7,17\}$\\
 &$n\equiv0(\mod2)$, $n\ge 4$, $(n,q)\neq(4,2)$ & $4$ & $\{p,r_{n-1},r_{2n-2}, r_{2n}\}$\\
 &$n\equiv1(\mod2)$, $n> 4$& $3$ &$\{p,r_{2n-2},r_{2n}\}$\\
 \hline
\end{tabular}}
\end{tab}
\vfill

\newpage
\begin{tab}\label{pindependenceex}{\bfseries $p$-independence numbers for finite
simple exceptional groups of Lie type}\vspace{1\baselineskip}

{\small \begin{tabular}{|c|l|c|c|}\hline
$G$&Conditions&$t(p,G)$&$\rho(p,G)$\\ \hline
$G_2(q)$&$q>2$&3&$\{p,r_3,r_6\}$\\ \hline
$F_4(q)$&none&3&$\{p,r_8,r_{12}\}$\\ \hline
$E_6(q)$&none&4&$\{p,r_8,r_9,r_{12}\}$\\ \hline
${^2E_6(q)}$&none&4&$\{p,r_8,r_{12},r_{18}\}$\\ \hline
$E_7(q)$&none&5&$\{p,r_7,r_9,r_{14},r_{18}\}$\\ \hline
$E_8(q)$&none&5&$\{p,r_{15},r_{20},r_{24},r_{30}\}$\\ \hline
$^3D_4(q)$&none&2&$\{p,r_{12}\}$\\\hline
${^2B_2(2^{2n+1})}$&$n\ge1$&4&$\{2,s_1,s_2,s_3\}$, where\\ &&&$s_1$
divides $2^{2n+1}-1$, \\ &&&$s_2$ divides $2^{2n+1}-2^{n+1}+1$\\
&&&$s_3$ divides $2^{2n+1}+2^{n+1}+1$\\ \hline
${^2G_2(3^{2n+1})}$&$n\ge1$&5&$\{3,s_1,s_2,s_3,s_4\}$, where\\
&&&$s_1\not=2$ divides $3^{2n+1}-1$\\ &&&$s_2\not=2$ divides
$3^{2n+1}+1$\\ &&&$s_3$ divides $3^{2n+1}-3^{n+1}+1$\\ &&&$s_4$
divides $3^{2n+1}+3^{n+1}+1$\\ \hline
${^2F_4(2^{2n+1})}$&$n\ge1$&4&$\{2,s_1,s_2,s_3\}$, where\\ &&&$s_1\not=3$ and
divides $2^{4n+2}-2^{2n+1}+1$\\ &&&$s_2$ divides
$2^{4n+2}+2^{3n+2}+2^{2n+1}+2^{n+1}+1$\\ &&&$s_3$ divides
$2^{4n+2}-2^{3n+2}+2^{2n+1}-2^{n+1}+1$\\ \hline ${}^2F_4(2)'$ &none&2 &
$\{2,13\}$\\ \hline
\end{tabular}}
\end{tab}

\newpage
\vfill
\begin{tab}\label{2independenceclass}{\bfseries 2-independence numbers for finite
simple classical groups of characteristic~${p\not=2}$}\vspace{1\baselineskip}

{\small
\begin{tabular}{|c|l|c|c|}
 \hline
 $G$&Conditions&$t(2,G)$&$\rho(2,G)$\\ \hline
 $A_{n-1}(q)$& $n=2$, $q\equiv1(\mod4)$& $3$ &$\{2,r_2, p\}$\\
 & $n=2$, $q\equiv3(\mod4)$, $q\neq 3$ & $3$ &$\{2,r_1, p\}$\\
  & $n\ge3$ and $n_2<(q-1)_2$& $2$ &$\{2,r_n\}$\\
 & $n\ge3$ and either $n_2>(q-1)_2$,& $2$ &$\{2,r_{n-1}\}$\\
 & or $n_2=(q-1)_2=2$&&\\
 & $2<n_2=(q-1)_2$& $3$ &$\{2,r_{n-1},r_n\}$\\
 \hline
 ${}^2A_{n-1}(q)$& $n_2>(q+1)_2$ & $2$& $\{2,r_{2n-2}\}$ \\
 & $n_2=1$ & $2$& $\{2,r_{2n}\}$ \\
  & $2<n_2<(q+1)_2$ & $2$& $\{2,r_{n}\}$ \\
 & $n\ge 3$, $2=n_2\le(q+1)_2$ & $2$& $\{2,r_{n/2}\}$ \\
 & $2<n_2=(q+1)_2$& $3$ & $\{2,r_{2n-2}, r_n\}$ \\
 \hline
 $B_n(q)$ or& $n>1$ is odd and $(q-1)_2=2$& $2$ & $\{2,r_n\}$\\
 $C_n(q)$& $n$ is even or $(q-1)_2>2$& $2$ &$\{2,r_{2n}\}$\\
 \hline
 $D_n(q)$ &$n\equiv0(\mod2)$, $n\ge 4$, $q\equiv 3(\mod 4)$ & $2$ &$\{2,r_{n-1}\}$\\
 &$n\equiv0(\mod2)$, $n\ge 4$, $q\equiv 1(\mod 4)$ & $2$ &$\{2,r_{2n-2}\}$\\
 &$n\equiv1(\mod2)$, $n>4$, $q\equiv 3(\mod 4)$ & $2$ & $\{2,r_n\}$\\
 &$n\equiv1(\mod2)$, $n>4$, $q\equiv 1(\mod 8)$& $2$ & $\{2,r_{2n-2}\}$\\
 &$n\equiv1(\mod2)$, $n>4$, $q\equiv 5(\mod 8)$& $3$ &$\{2,r_n,r_{2n-2}\}$\\
 \hline
 ${^2D_n(q)}$&$n\equiv0(\mod2)$, $n\ge 4$& $2$ & $\{2,r_{2n}\}$\\
  &$n\equiv1(\mod2)$, $n>4$, $q\equiv 1(\mod 4)$ & $2$ & $\{2,r_{2n}\}$\\
 &$n\equiv1(\mod2)$, $n>4$, $q\equiv 7(\mod 8)$& $2$ & $\{2,r_{2n-2}\}$\\
 &$n\equiv1(\mod2)$, $n>4$, $q\equiv 3(\mod 8)$& $3$ & $\{2,r_{2n-2},r_{2n}\}$\\
 \hline
\end{tabular}}
\end{tab}

\vfill

\begin{tab}\label{2independenceex}{\bfseries 2-independence numbers for finite
simple exceptional groups of Lie type
of characteristic~${p\not=2}$}\vspace{1\baselineskip}

{\small \begin{tabular}{|c|l|c|c|}
 \hline
 $G$&Conditions&$t(2,G)$&$\rho(2,G)$\\ \hline
$G_2(q)$&none&3&$\{2,r_3,r_6\}$\\ \hline
$F_4(q)$&none&2&$\{2,r_{12}\}$\\ \hline
$E_6(q)$&none&3&$\{2,r_9,r_{12}\}$\\ \hline
${^2E_6(q)}$&none&3&$\{2,r_{12},r_{18}\}$\\ \hline
$E_7(q)$&$q\equiv1(\mod4)$&3&$\{2,r_{14},r_{18}\}$\\
&$q\equiv3(\mod4)$&3&$\{2,r_7,r_9\}$\\ \hline $E_8(q)$&none&5&$\{2,r_{15},r_{20},r_{24},r_{30}\}$\\
\hline ${^3D_4(q)}$&none&2&$\{2,r_{12}\}$\\ \hline
${^2G_2(3^{2n+1})}$&none&3&$\{2,s_1,s_2\}$, where\\
&&&$s_1$ divides $3^{2n+1}-3^{n+1}+1$,\\
&&&$s_2$ divides $3^{2n+1}+3^{n+1}+1$.\\  \hline
\end{tabular}}
\end{tab}
\vfill
\newpage

\begin{tab}\label{independence}{\bfseries Independence numbers for finite
     simple classical groups}

{\small    \begin{tabular}{|c|l|c|c|}
  \hline
  $G$ & Conditions & $t(G)$ & $\rho(G)$\\
  \hline
 $A_{n-1}(q)$ & $n=2$, $q>3$ & $3$ & $\{p,r_1,r_2\}$\\
 & $n=3$, $(q-1)_3=3\text{ and }q+1\neq 2^k$ & $4$ & $\{p,3,r_2,r_3\}$\\
 & $n=3$, $(q-1)_3\neq 3\text{ and }q+1\neq 2^k$ & $3$ & $\{p,r_2,r_3\}$\\
 & $n=3$, $(q-1)_3=3\text{ and }q+1=2^k$ & $3$ & $\{p,3,r_3\}$\\
 & $n=3$, $(q-1)_3\neq 3\text{ and }q+1=2^k$ & $2$ & $\{p,r_3\}$\\
 & $n=4$ & $3$ & $\{p,r_{n-1},r_n\}$\\
 & $n=5,6$, $q=2$ & $3$ & $\{5,7,31\}$\\
 & $7\le n\le 11$, $q=2$ & $\left[\frac{n-1}{2}\right]$ & $\{r_i~|~i\neq 6, \left[\frac{n}{2}\right]<i\le n\}$\\
 & $n\ge 5$ and $q>2$ or  $n\ge 12$ and $q=2$& $\left[\frac{n+1}{2}\right]$ & $\{r_i~|~\left[\frac{n}{2}\right]<i\le n\}$\\
  \hline
 ${}^2A_{n-1}(q)$ & $n=3,q\not=2, (q+1)_3=3,\text{ and }q-1\neq 2^k$ & $4$ & $\{p,3,r_1,r_6\}$\\
 & $n=3$, $(q+1)_3\neq 3\text{ and }q-1\neq 2^k$ & $3$ & $\{p,r_1,r_6\}$\\
 & $n=3$, $(q+1)_3=3\text{ and }q-1=2^k$ & $3$ & $\{p,3,r_6\}$\\
 & $n=3$, $(q+1)_3\neq 3\text{ and }q-1=2^k$ & $2$ & $\{p,r_6\}$\\
 & $n=4$, $q=2$ & $2$ & $\{2,5\}$\\
 & $n=4$, $q>2$ & $3 $ & $\{p,r_4,r_6\}$\\
 & $n=5$, $q=2$ & $3$ & $\{2,5,11\}$\\
 & $n\ge 5$ and $(n,q)\neq (5,2)$ & $\left[\frac{n+1}{2}\right]$  & $\{r_{i/2}~|~\left[\frac{n}{2}\right]<i\le n,$\\
 & & & $i\equiv{2}(\mod 4)\}\cup$\\
 & & & $\cup\{r_{2i}~\vert~\left[\frac{n}{2}\right]<i\le n,$\\
 & & & $i\equiv{1}(\mod 2)\}$\\
 & & & $\cup\{r_i~|~\left[\frac{n}{2}\right]< i\le n,$\\
 & & & $i\equiv{0}(\mod 4)\}$\\
 \hline
 $B_n(q)$ or & $n=2$, $q>2$ & $2$ & $\{p, r_4\}$\\
 $C_n(q)$ & $n=3$ and $q=2$ & $2$ & $\{5,7\}$\\
 & $n=4$ and $q=2$ & $3$ & $\{5,7,17\}$\\
 & $n=5$ and $q=2$ & $4$ & $\{7,11,17,31\}$\\
 & $n=6$ and $q=2$ & $5$ & $\{7,11,13,17,31\}$\\
 & $n>2$,  $(n,q)\neq(3,2),(4,2),(5,2),(6,2)$ & $\left[\frac{3n+5}{4}\right]$ & $\{r_{2i}~|~\left[\frac{n+1}{2}\right]\le i\le n\}\cup$\\
 & &  &$\cup\{r_{i}~|~\left[\frac{n}{2}\right]<i\le n,$ \\
 & & & $i\equiv{1}(\mod 2)\}$\\ \hline
 $D_n(q)$ & $n=4$ and $q=2$ & $2$ & $\{5,7\}$\\
 & $n=5$ and $q=2$ & $4$ & $\{5,7,17,31\}$\\
 & $n=6$ and $q=2$ & $4$ & $\{7,11,17,31\}$\\
 &$n\ge 4$, $n\not\equiv3\pmod4$,&$\left[\frac{3n+1}{4}\right]$&$\{r_{2i}\mid
\left[\frac{n+1}{2}\right]\le
i< n\}\cup$\\
&$(n,q)\not=(4,2),(5,2),(6,2)$&&$\cup\{r_i\mid \left[\frac{n}{2}\right]<i\le
n,$\\
&&&$i\equiv 1\pmod2\}$\\
&$n\equiv3\pmod4$&$\frac{3n+3}{4}$&$\{r_{2i}\mid \left[\frac{n+1}{2}\right]\le
i< n\}\cup$\\
&&&$\cup\{r_i\mid \left[\frac{n}{2}\right]\le i\le
n,$\\ \hline
 ${}^2D_n(q)$ & $n=4$ and $q=2$ & $3$ & $\{5,7,17\}$\\
 & $n=5$ and $q=2$ & $3$ & $\{7,11,17\}$\\
 & $n=6$ and $q=2$ & $5$ & $\{7,11,13,17,31\}$\\
 & $n=7$ and $q=2$ & $5$ & $\{11,13,17,31,43\}$\\
&$n\ge 4$, $n\not\equiv1\pmod
4,$&$\left[\frac{3n+4}{4}\right]$&$\{r_{2i}\mid \left[\frac{n}{2}\right]\le
i\le n\}\cup$\\
&$(n,q)\not=(4,2),(6,2),(7,2)$&&$\cup\{r_i\mid \left[\frac{n}{2}\right]<i<
n,$\\
&&&$i\equiv 1\pmod2\}$\\
 & $n>4$, $n\equiv 1(\mod 4)$, $(n,q)\neq(5,2)$ & $\left[\frac{3n+4}{4}\right]$ & $\{r_{2i}~|~\left[\frac{n}{2}\right]< i\le n\}\cup$\\
 & &  &$\cup\{r_{i}~|~\left[\frac{n}{2}\right]<i< n,$ \\
 & & & $i\equiv{1}(\mod 2)\}$\\
 \hline
\end{tabular}}
\end{tab}

\newpage

\begin{tab}\label{indepexcept}{\bfseries Independence numbers for
finite simple exceptional groups of Lie type}\vspace{1\baselineskip}

{\small\begin{tabular}{|c|l|c|c|}
  \hline
  $G$ & Conditions &  $t(G)$&$\rho(G)$\\
  \hline
  $G_2(q)$ & $q>2$&3&  $\{p,r_3,r_6\}$ \\ \hline
$F_4(q)$&$q=2$&4&$\{5,7,13,17\}$\\
&$q>2$&5&$\{r_3,r_4,r_6,r_8,r_{12}\}$\\ \hline
$E_6(q)$& &5&$\{r_4,r_5,r_8,r_9,r_{12}\}$\\
\hline
${}^2E_6(q)$&&5&$\{r_4,r_8,r_{10},r_{12},r_{18}\}$\\ \hline
$E_7(q)$&&8&$\{r_5,r_7,r_8,r_9,r_{10},r_{12},r_{14},r_{18}\}$\\ \hline
$E_8(q)$&&12&$\{r_5,r_7,r_8,r_9,r_{10},r_{12},r_{14},r_{15},r_{18},r_{20},r_{
24 } ,
r_{30}\} $\\ \hline
${}^3D_4(q)$&$q=2$&2&$\{2,13\}$\\
&$q>2$&3&$\{r_3,r_6,r_{12}\}$\\ \hline
${^2B_2(2^{2n+1})}$&$n\ge1$&4&$\{2,s_1,s_2,s_3\}$,
where\\
&&&$s_1$ divides $2^{2n+1}-1$, \\
&&&$s_2$ divides $2^{2n+1}-2^{n+1}+1$\\ &&&$s_3$ divides $2^{2n+1}+2^{n+1}+1$\\
\hline
${^2G_2(3^{2n+1})}$&$n\ge1$&5&$\{3,s_1,s_2,s_3,s_4\}$, where\\ &&&$s_1\not=2$
divides $3^{2n+1}-1$\\
&&&$s_2\not=2$ divides $3^{2n+1}+1$\\ &&&$s_3$ divides $3^{2n+1}-3^{n+1}+1$\\
&&&$s_4$ divides $3^{2n+1}+3^{n+1}+1$\\ \hline
${}^2F_4(2^{2n+1})$&$n\ge2$,&$5$&$\{s_1,s_2,s_3,s_4,s_5\}$,
  where\\
&&&$s_1$ and divides $2^{2n+1}+1$,\\
&&&$s_2$ divides $2^{4n+2}+1$,\\
&&&$s_3\not=3$ and divides $2^{4n+2}-2^{2n+1}+1$,\\
&&&$s_4$ divides $2^{4n+2}-2^{3n+2}+2^{2n+1}-2^{n+1}+1$,\\
&&&$s_5$ divides $2^{4n+2}+2^{3n+2}+2^{2n+1}+2^{n+1}+1$,\\ \hline
  ${}^2F_4(2)'$ &none&  $3$&$\{3,5,13\}$\\ \hline
 ${}^2F_4(8)$&none&4&$\{7,19,37,109\}$ \\ \hline
\end{tabular}}
\end{tab}

\end{document}